\documentclass{fic-l}

\newcommand{\C}{{\mathbb C}}
\newcommand{\R}{{\mathbb R}}
\newcommand{\Z}{{\mathbb Z}}
\newcommand{\Q}{{\mathbb Q}}
\renewcommand{\P}{{\mathbb P}}
\newcommand{\s}{{\mathbb S}}
\newcommand{\B}{{\mathbb B}}
\newcommand{\I}{{\mathbb I}}

              \newcommand{\J}{{\bf J}}
              \newcommand{\M}{{\bf M}}
              
              \newcommand{\U}{{\bf U}}
              \newcommand{\T}{{\bf T}}
              
              \newcommand{\E}{{\bf E}}
              
              \renewcommand{\L}{{\cal L}}
              \renewcommand{\O}{{\cal O}}
              
              \newcommand{\G}{{\bf G}}

\newtheorem{theorem}{Theorem}[section]
\newtheorem{lemma}[theorem]{Lemma}

\newtheorem{problem}[theorem]{Problem}
\newtheorem{conjecture}[theorem]{Conjecture}
\newtheorem{definition}[theorem]{Definition}

\theoremstyle{definition}
\newtheorem{example}[theorem]{Example}

\theoremstyle{remark}
\newtheorem{remark}[theorem]{Remark}

\numberwithin{equation}{section}

\usepackage{amscd}

\begin{document}

\title{Pseudoholomorphic curves in four-orbifolds and some applications}

\author{Weimin Chen}
\address{Mathematics Department\\Tulane University\\
New Orleans, LA 70118, USA\\}
\email{wchen@math.tulane.edu}
\thanks{The author was supported in part by NSF grant
DMS-0304956.}

\date{December 1, 2004}



\maketitle

The main purpose of this paper is to summarize the basic ingredients,
illustrated with examples, of a pseudoholomorphic curve theory for symplectic
$4$-orbifolds. These are extensions of relevant work of Gromov, McDuff
and Taubes on symplectic $4$-manifolds concerning pseudoholomorphic curves
and Seiberg-Witten theory (cf. \cite{Gr, McD1, McD2, T1, T3}). They
form the technical backbone of \cite{C2, C3} where it was shown that a
symplectic $s$-cobordism of elliptic $3$-manifolds (with a canonical contact
structure on the boundary) is smoothly a product. We believe that this theory
has a broader interest and may find applications in other problems. One
interesting feature is that existence of pseudoholomorphic curves gives
certain restrictions on the singular points of the $4$-orbifold
contained by the pseudoholomorphic curves.

There are four sections. The first one is concerned with the Fredholm
theory for pseudoholomorphic curves in symplectic orbifolds, which is based
on the theory of maps of orbifolds developed in \cite{C1}, particularly the
topological structure of the corresponding mapping spaces. In the second
section, we discuss the orbifold version of adjunction formula and a formula
expressing the intersection number of two distinct pseudoholomorphic curves
in terms of contributions from each point in the intersection. Unlike the
manifold case, the contribution from a singular point may not be integral,
but it can be determined from the local structure of the pseudoholomorphic
curves near the singular point, which is important in applications. The third
section is occupied by a brief summary of the Seiberg-Witten theory for smooth
$4$-orbifolds and a discussion on the orbifold analog of the theorems of
Taubes concerning existence of pseudoholomorphic curves. We also
include here an example to illustrate how a combination of the various
ingredients is put into use. In the last section, we discuss some
applications and work in progress.

\section{Moduli space of pseudoholomorphic maps}
We begin with a brief review on the theory of maps of orbifolds
developed in \cite{C1}. (For an earlier version of this work, see
\cite{C0}.)

First of all, we recall the definition of (smooth) orbifolds from Satake
\cite{Sa} (see also \cite{Ka0}). Let $X$ be a paracompact, Hausdorff space.
An $n$-dimensional orbifold structure on $X$ is given by an atlas of local
charts $\U$, where $\U=\{U_i|i\in I\}$ is an open cover of $X$ satisfying
the following condition.
$$
(\ast) \mbox{    For any }p\in U_i\cap U_j, U_i,U_i\in\U,
\mbox{ there is a }
U_k\in\U \mbox{ such that }p\in U_k\subset U_i\cap U_j.
$$
Moreover, (1) for each $U_i\in\U$, there exists a triple $(V_i,G_i,\pi_i)$,
called a local uniformizing system, where $V_i$ is an $n$-dimensional smooth
manifold, $G_i$ is a finite group acting smoothly and effectively on $V_i$, and
$\pi_i:V_i\rightarrow U_i$ is a continuous map inducing a homeomorphism
$U_i\cong V_i/G_i$, (2) for each pair $U_i,U_j\in\U$ with $U_i\subset U_j$,
there is a set $Inj(U_i,U_j)=\{(\phi,\lambda)\}$, whose elements are called
injections, where $\phi:V_i\rightarrow V_j$ is a smooth open embedding and
$\lambda:G_i\rightarrow G_j$ is an injective homomorphism,
such that $\phi$ is $\lambda$-equivariant and satisfies $\pi_i=\pi_j\circ\phi$,
and $G_i\times G_j$ acts transitively on $Inj(U_i,U_j)$ by
$$
(g,g^\prime)\cdot (\phi,\lambda)=(g^\prime\circ\phi\circ g^{-1},
Ad(g^\prime)\circ\lambda\circ Ad(g^{-1})), \forall g\in G_i,
g^\prime\in G_j, (\phi,\lambda)\in Inj(U_i,U_j),
$$
(3) the injections are closed under composition for any $U_i,U_j,U_k\in\U$
with $U_i\subset U_j\subset U_k$.

Notice that for any refinement of $\U$ satisfying $(\ast)$, there is an induced
orbifold structure on $X$. Two orbifold structures on $X$ are called equivalent
if they induce isomorphic orbifold structures. With the preceding understood,
an orbifold is a paracompact, Hausdorff space with an equivalence class of
orbifold structures. One can similarly define orbifolds with boundary,
fiber bundles over orbifolds, etc. See \cite{Sa} for more details.

Let $X$ be an orbifold, and $\U$ be an orbifold structure on $X$. One may
assume without loss of generality that for any $p\in X$, there is a
$U_p\in\U$ centered
at $p$ in the sense that in the corresponding local uniformizing system
$(V_p,G_p,\pi_p)$, $V_p\subset\R^n$ is an open ball centered at $0$ with
$\pi_p(0)=p$ and $G_p$ acts linearly on $V_p$. The group $G_p$ is called the
isotropy group at $p$, and $p\in X$ is called a singular (resp. smooth) point
if $G_p\neq \{e\}$ (resp. $G_p=\{e\}$).

\begin{example}
A natural class of orbifolds is provided by the global quotients $X=Y/G$,
where $Y$ is a smooth manifold and $G$ is a discrete group acting smoothly
and effectively on $Y$ with only finite isotropy. Such orbifolds are 
called ``good'' orbifolds in Thurston \cite{Th}.

For example, consider the orbifold $X=\s^2/\Z_n$, where $\s^2\subset\R^3$
is the unit sphere, and the $\Z_n$-action is generated
by a rotation of angle $2\pi/n$ about the $z$-axis. There are two singular
points in $X$, the ``north pole'' and the ``south pole'', which are the
orbits of $(0,0,1)$ and $(0,0,-1)$ under the $\Z_n$-action respectively,
both having isotropy group $\Z_n$.

Not every orbifold is good. For example, consider the ``tear-drop''
orbifold in Thurston \cite{Th}, where the space $X=\s^2$, and the orbifold has
only one singular point $p$ such that a local uniformizing system centered
at $p$ is given by $(D,\Z_n,\pi)$, where $D\subset\C$ is a disc centered
at $0$, the action of $\Z_n$ is generated by a rotation of angle $2\pi/n$,
and $\pi(z)=z^n, \forall z\in D$.

Both orbifolds, $\s^2/\Z_n$ and the ``tear-drop'' orbifold, are examples of
orbifold Riemann surfaces. A general orbifold Riemann surface is
an orbifold where the underlying space is a Riemann surface and the orbifold
has finitely many singular points, each having a local uniformizing system
defined as in the case of the ``tear-drop'' orbifold above.
\hfill $\Box$

\end{example}

The concept of orbifold was introduced by Satake \cite{Sa} under
the name ``V-manifold'', where orbifolds were perceived as a class
of singular spaces to which the usual differential geometry on smooth
manifolds can be suitably extended by working locally and equivariantly
on each local uniformizing system. For example, a differential form on an
orbifold is a family of equivariant differential forms on the local
uniformizing systems which are compatible under the injections. The concept
was later rediscovered by Thurston \cite{Th}, who also invented the
more popular name ``orbifold''. Thurston was more concerned with the
topological aspects of orbifolds. For instance, Thurston introduced through a
novel covering theory the important concept of orbifold fundamental
group, which is different from the usual fundamental group. For example,
the orbifold fundamental groups of $\s^2/\Z_n$ and the ``tear-drop'' orbifold
in Example 1.1 are $\Z_n$ and the trivial group respectively, even though
both orbifolds have the same underlying space $\s^2$.
Orbifolds can be defined equivalently using the categorical
language of \'{e}tale topological groupoids, where the orbifold
fundamental group may be interpreted as the fundamental group of
an \'{e}tale topological groupoid (cf. eg. \cite{Hae, BH}).

Despite the fact that the concept of orbifold was extensively used
in the literature, there has not been a well-developed, satisfactory theory
of maps for orbifolds. The notion of maps that has been used mostly in the
literature, which is usually called a $C^\infty$ map (sometime it is
called an orbifold map), is essentially the notion of V-manifold map
introduced by Satake in \cite{Sa}, which is roughly speaking a continuous map
that can be locally lifted to a smooth map between local uniformizing systems.
Although such a notion of maps is sufficient for the purpose of defining
diffeomorphisms of orbifolds, or defining sections of an orbifold vector
bundle, etc., cf. Remark 1.7 (3) below, it has serious drawbacks in the fact
that the space of maps does
not have a well-understood topological structure in general, and that it is
not compatible with the distinct topological structure of orbifolds discovered
by Thurston which we alluded to earlier. On the other hand, the topological
aspect of orbifolds has been dealt with mainly through the classifying spaces
of topological groupoids (cf. eg. \cite{Hae}), rather than building upon a
satisfactory theory of maps as we do in the manifold case.

With the preceding understood, we now turn to the theory of maps
of orbifolds developed in \cite{C1}. The main advantage of
this theory lies in the fact that it accommodates at the same time both
the differential geometric and the topological aspects of orbifolds.
More concretely, in Part I of \cite{C1}, it was shown that the space
of maps in this theory has a natural infinite dimensional orbifold
structure, hence it is amenable to techniques of differential
geometry and global analysis on orbifolds. On the other hand, in
Part II of \cite{C1}, a basic homotopy theory and an analog of CW-complex
theory were developed (based upon the use of maps rather than classifying
spaces), which provide a machinery for studying homotopy class of maps
between orbifolds, and substantially extend the existing homotopy theory
defined via classifying spaces of topological groupoids (cf. \cite{Hae},
compare also \cite{Moer}).

The theory of maps in \cite{C1} was developed in the more general
context of orbispaces which are generalizations of orbifolds. For
technical reason we adopted a modified formalism, cf. Remark 1.7
(4) below.

\begin{definition} {\em(}See \cite{C1}, also \cite{C0}{\em)} 
Let $X$ be a locally connected topological space. An `orbispace
structure' on $X$ is given by an atlas of local charts $\U$, where
$\U=\{U_i|i\in I\}$ is an open cover of $X$ which satisfies the
following condition.
\begin{itemize}
\item [{(\#)}] Each $U_i\in\U$ is connected, and for any $U_i\in\U$, a
connected open subset of $U_i$ is also an element of $\U$.
\end{itemize}
Moreover,
\begin{itemize}
\item for each $U_i\in\U$, there is a triple $(\widehat{U_i},G_{U_i},
\pi_{U_i})$, where $\widehat{U_i}$ is a locally connected and connected
topological space with a continuous left action of a discrete
group $G_{U_i}$, and $\pi_{U_i}:\widehat{U_i}\rightarrow U_i$ is
a continuous map inducing a homeomorphism $\widehat{U_i}/G_{U_i}\cong
U_i$,
\item there exists a family of discrete sets
$$
\T=\{T(U_i,U_j)|U_i,U_j\in\U \mbox{ such that } U_i\cap
U_j\neq\emptyset\}
$$
which obeys
\begin{itemize}
\item $T(U_i,U_i)=G_{U_i}$, and for any $i\neq j$, if we let
$\{W_u|u\in I_{ij}\}$ be the set of connected components of 
$U_i\cap U_j$, then there are discrete sets $T_{W_u}(U_i,U_j)$ 
such that $T(U_i,U_j)=\sqcup_{u\in I_{ij}}T_{W_u}(U_i,U_j)$.
\item Each $\xi\in T_{W_u}(U_i,U_j)$ is assigned with a
homeomorphism $\phi_\xi$, whose domain and range are connected
components of the inverse image of $W_u$ in $\widehat{U_i}$ and
$\widehat{U_j}$ respectively, such that
$\pi_{U_i}|_{\mbox{Domain }(\phi_\xi)}=\pi_{U_j}\circ\phi_\xi$.
For each $\xi\in T(U_i,U_i)=G_{U_i}$, $\phi_\xi$ is the
self-homeomorphism of $\widehat{U_i}$ induced by the action of
$\xi$.
\item For any $\xi\in T(U_i,U_j)$, $\eta\in T(U_j,U_k)$, and any $x\in
\mbox{Domain }(\phi_\xi)$ such that $\phi_\xi(x)\in\mbox{Domain
}(\phi_\eta)$, there exists an $\eta\circ\xi(x)\in T(U_i,U_k)$
such that $\phi_{\eta\circ\xi(x)}(x)=\phi_{\eta}(\phi_\xi(x))$.
Moreover, $x\mapsto \eta\circ\xi(x)$ is locally constant, and the
composition $(\xi,\eta)\mapsto\eta\circ\xi(x)$ is associative and
coincides with the group multiplication in $G_{U_i}$ when
restricted to $T(U_i,U_i)=G_{U_i}$.
\item Each $\xi\in T(U_i,U_j)$ has an inverse $\xi^{-1}\in
T(U_j,U_i)$, with $\mbox{Domain }(\phi_\xi)
=\mbox{Range}(\phi_{\xi^{-1}})$, $\mbox{Domain }(\phi_{\xi^{-1}})
=\mbox{Range }(\phi_\xi)$, and $\xi^{-1}\circ\xi(x)=e\in G_{U_i}
\forall x\in\mbox{Domain }(\phi_\xi)$, $\xi\circ\xi^{-1}(x)=e
\in G_{U_j}\forall x\in \mbox{Domain }(\phi_{\xi^{-1}})$.
\end{itemize}
\end{itemize}

For any refinement of $\U$ satisfying $(\#)$, there is an induced
orbispace structure on $X$. Two orbispace structures on $X$ are
called `equivalent' if they induce isomorphic orbispace structures.
Finally, an `orbispace' is a locally connected topological space
with an orbispace structure, which is uniquely determined up to
equivalence.
\end{definition}

\begin{remark}
(1) The orbifolds in \cite{C1} are defined in the formalism of
Definition 1.2. They differ from those defined in Satake \cite{Sa} in
two aspects. The first one which is more essential is that the group
action in a local uniformizing system is no longer required to be effective
in \cite{C1}. An orbifold in the classical sense (ie. the local
group actions are effective) is called `reduced' (cf. \cite{CR}).
The second aspect is the technical requirement that for any $U_i,U_j\in\U$ 
in the atlas of local charts, a connected component of $U_i\cap U_j$ is also 
an element of $\U$ (note that $T(U_i,U_j)$,
which corresponds to $Inj(U_i,U_j)$ in Satake's formalism when
$U_i\subset U_j$, is defined as long as $U_i\cap U_j\neq\emptyset$).
This technical assumption can always be met by replacing $\U$ with a
refinement of $\U$ which defines an equivalent orbifold structure (cf.
Proposition 2.1.3, Part I of \cite{C1}).

(2) The orbispaces in the sense of Definition 1.2 belong to a restricted
class of orbispaces in Haefliger \cite{Hae1}, which were defined using
the language of \'{e}tale topological groupoids. To be more precise,
if we interpret the orbispace structures in Definition 1.2 in terms
of \'{e}tale topological groupoids, the formalism of Definition 1.2 (which
is a modified version of Satake's) requires the corresponding \'{e}tale
topological groupoids to satisfy some additional conditions (ie. the
conditions (C1), (C2) in Part I of \cite{C1}). Moreover, simple examples show
that these conditions are not preserved under the Morita equivalence of
topological groupoids. Besides orbifolds, it is known that global quotient
orbispaces $X=Y/G$ (where $Y$ is a locally connected topological space 
and $G$ is a discrete group) and the orbihedra defined in \cite{Hae1}
satisfy the conditions of Definition 1.2. See \cite{C1} for more details.
\hfill $\Box$

\end{remark}

Let $X,X^\prime$ be any orbispaces in the sense of Definition
1.2, with orbispace structures $\U$, $\U^\prime$ respectively.
In order to define maps between $X$ and $X^\prime$, we consider sets of
data $(\{\hat{f}_\alpha\},\{\rho_{\beta\alpha}\})$ (which will be called
a `compatible system' below), where

\begin{itemize}
\item there is an open cover $\{U_\alpha\}$ of $X$ with $U_\alpha\in\U$,
and a correspondence $U_\alpha\mapsto U_\alpha^\prime\in\U^\prime$,
\item each $\hat{f}_\alpha:\widehat{U_\alpha}\rightarrow
\widehat{U_\alpha^\prime}$ is a continuous map, and each $\rho_{\beta\alpha}:
T(U_\alpha,U_\beta)\rightarrow T(U_\alpha^\prime,U_\beta^\prime)$ is
a mapping of discrete sets, such that
\begin{itemize}
\item [{(a)}] $\phi_{\rho_{\beta\alpha}(\xi)}\circ \hat{f}_\alpha(x)
=\hat{f}_\beta\circ\phi_\xi(x)$ for any $\alpha,\beta$, where $\xi\in
T(U_\alpha,U_\beta)$, $x\in\mbox{Domain }(\phi_\xi)$.
\item [{(b)}] $\rho_{\gamma\alpha}(\eta\circ\xi(x))=\rho_{\gamma\beta}(\eta)
\circ\rho_{\beta\alpha}(\xi)(\hat{f}_\alpha(x))$ for any
$\alpha,\beta,\gamma$, where $\xi\in T(U_\alpha,U_\beta)$,
$\eta\in T(U_\beta,U_\gamma)$, $x\in\phi_\xi^{-1}(\mbox{Domain }
(\phi_\eta))$. (Note: $\hat{f}_\alpha(x)\in\phi_{\rho_{\beta\alpha}(\xi)}^{-1}
(\mbox{Domain }(\phi_{\rho_{\gamma\beta}(\eta)})$, which follows from (a)
and the assumption that $x\in\phi_\xi^{-1}(\mbox{Domain }(\phi_\eta))$.)
\end{itemize}
(Observe that if set $\rho_\alpha\equiv\rho_{\alpha\alpha}:G_{U_\alpha}
\rightarrow G_{U_\alpha^\prime}$, then condition (b) above implies that
$\rho_\alpha$ is a homomorphism and condition (a) above implies that
$\hat{f}_\alpha$ is $\rho_\alpha$-equivariant.)
\end{itemize}

One can introduce an equivalence relation between compatible systems
as follows. Suppose $(\{\hat{f}_a\},\{\rho_{ba}\})$ is another compatible
system. We say that $(\{\hat{f}_a\},\{\rho_{ba}\})$ is `induced' from
$(\{\hat{f}_\alpha\},\{\rho_{\beta\alpha}\})$, if
there exists a set of data $(\theta,\{\xi_a\},\{\xi^\prime_a\})$, where
$\theta:a\mapsto\alpha$ is a mapping of indices satisfying $U_a\subset
U_{\theta(a)}$, and $\xi_a\in T(U_a,U_{\theta(a)})$, $\xi^\prime_a\in
T(U_a^\prime,U_{\theta(a)}^\prime)$, such that
\begin{itemize}
\item $\hat{f}_a=(\phi_{\xi_a^\prime})^{-1}\circ
\hat{f}_{\theta(a)}\circ\phi_{\xi_a}$, and
\item $\rho_{ba}(\eta)=(\xi_b^\prime)^{-1}\circ
\rho_{\theta(b)\theta(a)}(\theta(\eta))\circ\xi_a^\prime(x)$,
$\forall x\in \hat{f}_a(\mbox{Domain}(\phi_\eta))$, where
$\theta(\eta)\equiv\xi_b\circ\eta\circ\xi_a^{-1}(z)\in
T(U_{\theta(a)},U_{\theta(b)})$, $\forall z\in\phi_{\xi_a}(\mbox{Domain }
(\phi_{\eta}))$.
\end{itemize}

\begin{definition} {\em(}See \cite{C1}, also \cite{C0}{\em)} 
Two compatible systems are said to be `equivalent' if they induce a common
compatible system.
\end{definition}

A map between two orbispaces is defined to be an equivalence class of
compatible systems in the sense of Definition 1.4. (It is shown in
\cite{C1} that the relation in Definition 1.4 is indeed an equivalence
relation.) With such a notion of maps it is shown in
\cite{C1} that the set of orbispaces in the sense of Definition 1.2 form a
category. Observe that a map $f:X\rightarrow X^\prime$ between orbispaces
induces a continuous map between the underlying topological spaces,
which, when $X,X^\prime$ are reduced orbifolds and $f$ is smooth, is a
V-manifold map of Satake \cite{Sa}. We denote by $[X;X^\prime]$
the space of maps between $X$ and $X^\prime$. Notice that $[X;X^\prime]$
is uniquely determined up to a bijection by the equivalence class of
orbispace structures on $X,X^\prime$.

The central result in Part I of \cite{C1}, which is also the technical
foundation for the subsequent development in Part II of \cite{C1}, is
the following structural theorem.

\begin{theorem}
Let $X,X^\prime$ be any orbispaces with orbispace structures $\U,\U^\prime$
respectively. Then $[X;X^\prime]$ is naturally an orbispace in the sense of
Haefliger \cite{Hae1}, provided that the following conditions are satisfied
by $X$:
\begin{itemize}
\item [{(1)}] the underlying topological space is paracompact, locally
compact and Hausdorff,
\item [{(2)}] the orbispace structure $\U$ obeys: for each $U\in\U$, the space
$\widehat{U}$ is locally compact and Hausdorff, and $\pi_U:\widehat{U}
\rightarrow U$ is proper.
\end{itemize}
\end{theorem}

Specializing in the case of orbifolds, we have the following theorem.
Here $[X;X^\prime]^r$ denotes the space of $C^r$ maps between
smooth orbifolds $X$ and $X^\prime$, where a $C^r$ map is the equivalence
class of a compatible system $(\{\hat{f}_\alpha\},\{\rho_{\beta\alpha}\})$
such that each $\hat{f}_\alpha$ is a $C^r$ map between smooth manifolds.

\begin{theorem}
Let $X,X^\prime$ be smooth orbifolds where $X$ is compact and connected.
Then $[X;X^\prime]^r$ is naturally a Banach orbifold {\em(}Hausdorff and
second countable{\em)}, such that at any $f\in [X;X^\prime]^r$,
there is a uniformizing system $(V_f,G_f)$, where $V_f$ is an open ball
centered at $0$ in the space of $C^r$ sections of $f^\ast (TX^\prime)$,
and $G_f$ acts on $V_f$ linearly. Moreover, for any
$f\in [X;X^\prime]^r$ and any point $p$ in the image of $f$,
$G_f$ is naturally a subgroup of the isotropy group $G_p$ at $p$.
In particular, a $C^r$ map is a smooth point in $[X;X^\prime]^r$
if its image contains a smooth point of $X^\prime$.
\end{theorem}

\begin{remark}
(1) In Theorem 1.6, the orbifolds $X,X^\prime$ are not
necessarily reduced, and $X$ can be more generally an orbifold
with boundary. For any $C^r$ map $f:X\rightarrow X^\prime$
and orbifold vector bundle $E$ over $X^\prime$, the pull-back
bundle of $E$ via $f$, which is denoted by $f^\ast E$, is well-defined
as a $C^r$ orbifold vector bundle up to isomorphisms. Finally, by the image
of a map between orbifolds, we always mean the image of the induced
continuous map between the underlying spaces.

(2) A Banach orbifold is the straightforward generalization (ie., not in
the modified formalism of Definition 1.2) of Satake's orbifolds
to the infinite dimensional setting, except that the group action on each
local uniformizing system is not required to be effective. In particular,
the action of $G_f$ in Theorem 1.6 may not be effective in general.

(3) Suppose both $X,X^\prime$ are reduced. Let $f^\flat$ be
any V-manifold map from $X$ to $X^\prime$ such that the inverse image
of the smooth locus of $X^\prime$ is dense in $X$. Then $f^\flat$
uniquely determines a map $f:X\rightarrow X^\prime$ of orbifolds as follows.
Pick any set of local liftings $\{\hat{f}_\alpha\}$ of $f^\flat$,
there is a uniquely determined set of mappings $\{\rho_{\beta\alpha}\}$
such that $(\{\hat{f}_\alpha\},\{\rho_{\beta\alpha}\})$ form a
compatible system, whose equivalence class is easily seen to depend on
$f^\flat$ only (cf. \cite{CR}). This perhaps explains why the notion of
V-manifold map is sufficient as far as only diffeomorphisms of orbifolds or
sections of an orbifold vector bundle are concerned.

(4) The notion of maps in \cite{C1} has its root in
the construction of Gromov-Witten invariants of symplectic
orbifolds in \cite{CR}, where the main technical obstacle was the
fact that the notion of V-manifold map is not sufficient to define
pull-back bundles. The notion of compatible systems\footnote{We would
like to point out that the present definition of compatible systems,
which relies on the modified definition of orbifolds in Definition
1.2, is different from that in \cite{CR}, and this modified formalism
(along with Definition 1.4) is crucial in proving Theorems 1.5 and 1.6.}
was introduced in \cite{CR} as a device to define pull-back bundles,
and a V-manifold map was called a `good map' in \cite{CR} if it supports
such a compatible system. Moreover, two compatible systems in \cite{CR} were
defined to be `isomorphic' if they define isomorphic pull-back bundles.
For reduced orbifolds, one can show that this somewhat indirect,
implicit definition is actually equivalent to the direct, more explicitly
defined one in Definition 1.4, due to the fact explained in Remark 1.7 (3)
above. The bottom line is: The Gromov-Witten invariants based on
\cite{C1} are the same as those defined in \cite{CR}. However, the
treatment in \cite{CR} was {\it ad hoc} in nature and can be substantially
simplified in the framework of \cite{C1}.

(5) In terms of \'{e}tale topological groupoids, a compatible
system may be interpreted as a groupoid homomorphism, and a map
between orbispaces in \cite{C1} may be regarded as a groupoid homomorphism
modulo certain Morita equivalence of topological groupoids. There are several
related notions of maps in the literature, see Haefliger \cite{Hae0},
Hilsum-Skandalis \cite{HS}, Pronk \cite{Pr}, and Abramovich-Vistoli \cite{AV}
(in the context of Deligne-Mumford stacks). Compare also \cite{Moer}.
\hfill $\Box$
\end{remark}

We illustrate Theorem 1.6 with the following example.

\begin{example}

Let $X=Y/G$ and $X^\prime=Y^\prime/G^\prime$, where $Y$ is compact
and connected, and $G,G^\prime$ are finite. (Here the actions of
$G,G^\prime$ are not assumed to be effective). Consider
$$
[(Y,G);(Y^\prime,G^\prime)]^r \equiv \{(\hat{f},\rho)|
\hat{f}:Y\rightarrow Y^\prime \mbox{ a $\rho$-equiv. $C^r$ map, }
\rho:G\rightarrow G^\prime \mbox{ a homo.}\}
$$
which is a Banach manifold with a smooth action of $G^\prime$ defined by
$$
g\cdot (\hat{f},\rho)=(g\circ \hat{f}, Ad(g)\circ\rho),\;
\forall g\in G^\prime,
(\hat{f},\rho)\in [(Y,G);(Y^\prime,G^\prime)]^r.
$$
Note that each $(\hat{f},\rho)\in [(Y,G);(Y^\prime,G^\prime)]^r$
is naturally a compatible system. With this understood, the correspondence
$$
(\hat{f},\rho)\mapsto\mbox{ the equivalence class of } (\hat{f},\rho)
\mbox{ in } [X;X^\prime]^r
$$
identifies the Banach orbifold $[(Y,G);(Y^\prime,G^\prime)]^r/G^\prime$
as an open and closed subspace of $[X;X^\prime]^r$, which is a
proper subspace in general when $\pi_1(Y)$ is nontrivial. (In fact,
each $f\in [X;X^\prime]^r$ induces a homomorphism $f_\ast:\pi_1^{orb}(X)
\rightarrow\pi_1^{orb}(X^\prime)$ between orbifold fundamental
groups, and $f\in [(Y,G);(Y^\prime,G^\prime)]^r/G^\prime$ iff $f_\ast(\pi_1(Y))
\subset\pi_1(Y^\prime)$.)

A regular neighborhood of $(\hat{f},\rho)$ in $[(Y,G);(Y^\prime,G^\prime)]^r$
can be described as follows. Note that there is a $G$-bundle
$E_{(\hat{f},\rho)}
\rightarrow Y$ which is the pull-back of $TY^\prime$ by $(\hat{f},\rho)$. Let
$V_{(\hat{f},\rho)}$ be a sufficiently small open ball centered at the origin 
in
the space of $G$-sections of $E_{(\hat{f},\rho)}$ of $C^r$ class. The isotropy
subgroup $G_{(\hat{f},\rho)}$ at $(\hat{f},\rho)$, which consists of
$g\in G^\prime$ such that $g\cdot \hat{f}(p)=\hat{f}(p)$, $g\rho(h)=\rho(h)g$
for all $p\in Y$ and $h\in G$, acts linearly on $V_{(\hat{f},\rho)}$ by
$$
g\cdot s=g\circ s,\; \forall g\in G_{(\hat{f},\rho)}, s\in V_{(\hat{f},\rho)}.
$$
With this understood, a regular neighborhood of $(\hat{f},\rho)$ can be
identified with $V_{(\hat{f},\rho)}$ via the $G_{(\hat{f},\rho)}$-equivariant
map defined by
$$
s(p)\mapsto (\exp_{\hat{f}(p)}s(p),\rho), \; \mbox{ where }
p\in Y, s\in V_{(\hat{f},\rho)}.
$$
In other words, $(V_{(\hat{f},\rho)},G_{(\hat{f},\rho)})$ is a uniformizing
system at the image of $(\hat{f},\rho)$ in $[X;X^\prime]^r$.

Finally, observe that $G_{(\hat{f},\rho)}$ is a subgroup of the isotropy
subgroup at $\hat{f}(p)$ for any $p\in Y$.
\hfill $\Box$
\end{example}

With these preparations, we now discuss the Fredholm theory for
pseudoholomorphic curves in symplectic orbifolds. 


\begin{definition}

{\em(1)} A symplectic structure on an orbifold $X$ is a closed,
non-degenerate $2$-form $\omega$, which is given by a family of closed,
non-degenerate $2$-forms $\{\omega_i\}$ on $V_i$ for each local 
uniformizing system $(V_i,G_i,\pi_i)$, which are equivariant under the 
local group actions and compatible with respect to the injections.

{\em(2)} An almost complex structure on an orbifold $X$ is an endomorphism 
$J:TX\rightarrow TX$ with $J^2=-1$, which is given by a family of
endomorphisms $\{J_i:TV_i\rightarrow TV_i\}$ with $J_i^2=-1$ for each local 
uniformizing system $(V_i,G_i,\pi_i)$, which are equivariant under the 
local group actions and compatible with respect to the injections. Let
$\omega$ be any symplectic structure on $X$. An almost complex structure
$J$ is called $\omega$-compatible if $\omega(\cdot,J\cdot)$ defines
a Riemannian metric on $X$.

\end{definition}

(We remark that when $X=Y/G$ is a global quotient by a finite group, 
a symplectic (resp. almost complex) structure on the orbifold $X$ is
given by a $G$-equivariant symplectic (resp. almost complex) structure
on $Y$. )

\vspace{2mm}

As for the general setup, let $(X,\omega)$ be a closed, reduced symplectic 
orbifold and $J$ be a fixed $\omega$-compatible almost complex structure on 
$X$. We will be considering $J$-holomorphic maps (always assuming nonconstant,
ie. the image does not consist of a single point) from an orbifold Riemann
surface $(\Sigma,j)$ into $(X,J)$. (A map $f:\Sigma\rightarrow X$ is
called $J$-holomorphic if each $\hat{f}_\alpha$ in a corresponding
compatible system is $J$-holomorphic, ie., $d\hat{f}_\alpha+J_\alpha\circ
d\hat{f}_\alpha\circ j_\alpha=0$.) Here an additional assumption is made
that each homomorphism $\rho_\alpha\equiv\rho_{\alpha\alpha}$ in the
compatible systems of a $J$-holomorphic map is injective. However, we do
not assume that the orbifold Riemann surface $\Sigma$ is reduced.
Furthermore, we will adopt the following convention for $\Sigma$: a point
$z\in\Sigma$ where the isotropy group $G_z$ acts trivially is called a
`regular point'; otherwise, it is called an `orbifold point'.

The most relevant information concerning a $J$-holomorphic map is its
local structure at each point. More precisely, let $f:\Sigma\rightarrow X$
be a $J$-holomorphic map, let $z\in\Sigma$, $p\in X$ such that
$f(z)=p$, and let $(D_z,G_z)$, $(V_p,G_p)$ be local uniformizing systems
at $z$, $p$ respectively. (Note that when $\Sigma$ is not reduced, $G_z$
can be any finite group in general.) Then $f$ determines a germ of pairs
$(\hat{f}_z,\rho_z)$, where $\rho_z: G_z\rightarrow G_p$ is an injective
homomorphism, and $\hat{f}_z:D_z\rightarrow V_p$ is $\rho_z$-equivariant
and $J$-holomorphic. Moreover, $(\hat{f}_z,\rho_z)$ is unique up to
a change $(\hat{f}_z,\rho_z)\mapsto (g\circ \hat{f}_z,Ad(g)\circ\rho_z)$
for some $g\in G_p$. We will call $(\hat{f}_z,\rho_z)$ a `local
representative' of $f$ at $z$.

With Theorem 1.6, the Fredholm theory in the orbifold setting is in
complete analogy with that in the smooth setting. We shall illustrate
it by considering the case where $X$ has at most isolated singularities,
which suffices for most of the applications in dimension $4$. Notice that
in this case the orbifold Riemann surface $\Sigma$ is necessarily reduced
by our assumptions.

Fixing a sufficiently large $r>0$, we consider the space $[\Sigma;X]^r$,
which, under the assumption on $X$, may be regarded as a Banach manifold
because each nonconstant $J$-holomorphic map from $\Sigma$ to $X$ contains
a smooth point of $X$ in its image, hence is a smooth point in $[\Sigma;X]^r$
by Theorem 1.6. On the other hand, it is easily seen that the space $\J^r$
of $\omega$-compatible almost complex structures of $C^r$ class is always a
Banach manifold, with the tangent space at $J\in\J^r$ being the space of
$C^r$ sections $A$ of the orbifold vector bundle $\mbox{End }TX$,
which obeys (1) $AJ+JA=0$, (2) $A^t=A$. (Here the transpose $A^t$ is taken
with respect to $\omega(\cdot,J\cdot)$.)

We also need to consider the moduli space of complex structures on
the orbifold Riemann surface $\Sigma$, which can be identified with
the moduli space of complex structures on the corresponding marked
Riemann surface (with the orbifold points being the marked
points), cf. \cite{CR}. For simplicity, we assume $\Sigma$ is an
orbifold Riemann sphere of $k\leq 3$ orbifold points. Then $\Sigma$
has a unique complex structure $j$ with a $(6-2k)$-dimensional
automorphism group $\G$.

Now for any $(f,J)\in [\Sigma;X]^r\times \J^r$, let $\E_{(f,J)}$ be
the space of $C^{r-1}$ sections $s$ of the orbifold bundle $\mbox{Hom}
(T\Sigma,f^\ast(TX))$ such that $s\circ j=-J\circ s$. Then there
is a Banach bundle $\E\rightarrow [\Sigma;X]^r\times \J^r$ whose
fiber at $(f,J)$ is $\E_{(f,J)}$. The zero set of
$\underline{L}:[\Sigma;X]^r\times \J^r\rightarrow\E$, where
$$
\underline{L}(f,J)=df+J\circ df\circ j,
$$
consists of pairs $(f,J)$ such that $f$ is $J$-holomorphic. We set
$$
\widetilde{\M}\equiv \{(f,J)|(f,J)\in\underline{L}^{-1}(0), f
\mbox{ is not multiply covered}\},
$$
and for any $J\in\J^r$, set $\widetilde{\M}_J\equiv\{f|(f,J)
\in\widetilde{\M}\}$.
(Here we say $f$ is not multiply covered if the induced
map is not multiply covered.)
The following lemma follows by the same analysis as in the manifold
case (cf. eg. \cite{McDS}) with an application of the Sard-Smale
theorem \cite{Sm}.

\begin{lemma}
The subspace $\widetilde{\M}\subset [\Sigma;X]^r\times \J^r$ is a
Banach submanifold, and for a generic $J\in\J^r$ which may be chosen
smooth, $\widetilde{\M}_J$ is a smooth, finite dimensional manifold.
\end{lemma}

The dimension of $\widetilde{\M}_J$ can be calculated using Kawasaki's
orbifold index theorem \cite{Ka}, see also Lemma 3.2.4 in \cite{CR}.
We state the formula here for the case where $\dim X=4$
(and $\Sigma$ is a reduced orbifold Riemann sphere with $k\leq 3$ orbifold
points). Let $f\in\widetilde{\M}_J$ and let $z_i\in\Sigma$ be the
orbifold points with orders $m_i$. Suppose in a local
representative $(\hat{f}_{z_i},\rho_{z_i}):(D_i,\Z_{m_i})\rightarrow
(V_i,G_i)$ of $f$ at $z_i$, the action of $\rho_{z_i}(\mu_{m_i})$ on
$V_i$ (here $\mu_m\equiv\exp(\sqrt{-1}\frac{2\pi}{m})$) is given by
$$
\rho_{z_i}(\mu_{m_i})\cdot (w_1,w_2)=(\mu_{m_i}^{m_{i,1}}w_1,
\mu_{m_i}^{m_{i,2}}w_2), \; 0<m_{i,1},m_{i,2}<m_i.
$$
Then $\dim \widetilde{\M}_J=2\tilde{d}$ at $f$ where
$$
\tilde{d}=c_1(TX)\cdot f_\ast([\Sigma])+2-\sum_{i=1}^k
\frac{m_{i,1}+m_{i,2}}{m_i}\in \Z.
$$

The automorphism group $\G$ of the complex structure $j$ on $\Sigma$ acts
smoothly on $\widetilde{\M}_J$
(see the end of \S 3.3 in Part I of \cite{C1} for a general discussion),
which is free because the maps in $\widetilde{\M}_J$ are not
multiply covered. Consequently, $\M_J\equiv\widetilde{\M}_J/\G$ is also
a smooth manifold with $\dim\M_J=\dim\widetilde{\M}_J-(6-2k)$.

It is well-known that in the smooth case, Lemma 1.10 implies
nonexistence of certain $J$-holomorphic curves in a symplectic
$4$-manifold when $J$ is chosen generic. For instance, there exists
no embedded $J$-holomorphic $2$-sphere $C$ with $C^2=-2$ and $J$
generic. We will see that this is no longer true if there is a
finite group action and $J$ is chosen equivariant. However, existence of
$J$-holomorphic curves for a generic equivariant $J$ will put
certain restrictions on the fixed-point set of the action, as shown in the
following example.

\begin{example}

Let $\Z_n$ act on a symplectic $4$-manifold $Y$ preserving the symplectic
structure, and let $p_1,p_2\in Y$ be two isolated fixed points. Suppose for
any generic equivariant $J$, there exists a $\Z_n$-invariant, embedded
$J$-holomorphic $2$-sphere $C$ with $C^2=-2$ such that $p_1,p_2\in C$.
We apply Lemma 1.10 to the orbifold $X\equiv Y/\Z_n$ and $J$-holomorphic
maps $\Sigma\equiv\s^2/\Z_n\rightarrow X$. Then $\M_J$, which is a smooth
manifold by Lemma 1.10, must have $\dim\M_J=2d\geq 0$ at $C/\Z_n$. Here
by the aforementioned formula for $\dim\widetilde{\M}_J$,
$$
d=\frac{1}{n}c_1(TY)\cdot C+2-\frac{1+m_1}{n}-\frac{1+m_2}{n}-(3-2)
=1-\frac{1+m_1}{n}-\frac{1+m_2}{n},
$$
where $0<m_1,m_2<n$ are the weights of the complex representations
of $\Z_n$ at $p_1,p_2$ in the direction normal to $C$. Note that $d\in\Z$,
which implies that either $m_1+m_2+2=n$ or $m_1+m_2+2=2n$.
But the latter case, which implies $m_1=m_2=n-1$, can be ruled out
because $d\geq 0$.
\hfill $\Box$
\end{example}

Note that a discussion concerning Lemma 1.10 and Example 1.11 where
$\Sigma$ is a general orbifold Riemann surface may be found in
\cite{C4}.

Besides Lemma 1.10, there is another type of regularity criterion
which is the orbifold version of Lemma 3.3.3 in \cite{McDS}. See
Lemma 3.2 of \cite{C3}. Here the almost complex structure $J$ is not
necessarily generic.

\begin{lemma}
Let $f:\Sigma\rightarrow X$ be a $J$-holomorphic map from a reduced
orbifold Riemann surface into a symplectic $4$-orbifold with only isolated
singularities. Suppose for any $z\in\Sigma$, the map $\hat{f}_z$ in a
local representative $(\hat{f}_z,\rho_z)$ of $f$ at $z$ is embedded. Then
$f$ is a smooth point in the space of $J$-holomorphic maps if
$c_1(T\Sigma)\cdot [\Sigma]>0$ and $c_1(TX)\cdot f_\ast ([\Sigma])>0$.
\end{lemma}

Finally, we state an orbifold version of the Gromov compactness
theorem. For simplicity, we assume that the orbifold Riemann surface
$\Sigma$ is reduced.

\begin{theorem}
Let $(X,\omega)$ be any compact closed symplectic orbifold,
$J$ be any $\omega$-compatible almost complex structure, and
$0\neq A\in H_2(X;\R)$. For any sequence $f_n:\Sigma\rightarrow X$
of $J$-holomorphic maps with $(f_n)_\ast([\Sigma])=A$, there is a
subsequence, reparametrized if necessary and still denoted by $f_n$
for simplicity, having the following significance: There are at
most finitely many simple closed loops $\gamma_1,\cdots,\gamma_l
\subset\Sigma$ containing no orbifold points, and a nodal orbifold
Riemann surface $\Sigma^\prime=\cup_\nu\Sigma_\nu$ obtained by
collapsing $\gamma_1,\cdots,\gamma_l$, and a $J$-holomorphic map
$f:\Sigma^\prime\rightarrow X$ with $f_\ast([\Sigma])=A\in H_2(X;\R)$,
such that
\begin{itemize}
\item [{(1)}] $f_n$ converges to $f$ in $C^\infty$ on any given
compact subset in the complement of $\gamma_1,\cdots,\gamma_l$,
\item [{(2)}] if $z_\nu\in\Sigma_\nu,z_\omega\in\Sigma_\omega$ are
two distinct points with orders $m_\nu,m_\omega$ respectively, such that
$z_\nu,z_\omega$ are the image of the same simple
closed loop collapsed under $\Sigma\rightarrow\Sigma^\prime$, then
$m_\nu=m_\omega=m$, and there exist local representatives $(\hat{f}_{\nu},
\rho_{\nu})$, $(\hat{f}_{\omega},\rho_{\omega})$ of $f$ at
$z_\nu,z_\omega$, which obey $\rho_{\nu}(\mu_m)
=\rho_{\omega}(\mu_m)^{-1}$ {\em(}here $\mu_m\equiv
\exp(\sqrt{-1}\frac{2\pi}{m})${\em)},
\item [{(3)}] if $f$ is constant over $\Sigma_\nu$, then either the underlying
Riemann surface $|\Sigma_\nu|$ has nonzero genus, or $\Sigma_\nu$ contains at
least $3$ special points where a special point is either an orbifold point
inherited from $\Sigma$ or the image of some $\gamma_i$,
$i=1,\cdots,l$, under $\Sigma\rightarrow\Sigma^\prime$.
\end{itemize}
\end{theorem}

Theorem 1.13 can be proved in the framework of \cite{C1} along the lines
of Parker-Wolfson \cite{PW} or Ye \cite{Ye}, particularly,
with help of an orbifold version of Arzela-Ascoli Theorem. A version
of Theorem 1.13 was proved in \cite{CR} in an {\it ad hoc} manner
without using the theory in \cite{C1}.

\begin{theorem}[Orbifold Arzela-Ascoli Theorem, cf. \cite{C1}]
Let $X,X^\prime$ be any complete Riemannian orbifolds where $X$ is compact
{\em(}with or without boundary{\em)}. For any sequence of $C^r$ maps from $X$
to $X^\prime$ which have bounded $C^r$ norms, there is a subsequence which
converges to a $C^{r-1}$ map in the $C^{r-1}$ topology.
\end{theorem}

\section{Adjunction and intersection formulae}

In this section, we explain the adjunction and intersection formulae for
pseudoholomorphic curves in an almost complex $4$-orbifold, which were proved
in \cite{C2} and are the orbifold versions of the relevant results of
Gromov \cite{Gr} and McDuff \cite{McD1, McD2}.

First of all, to set the stage we let $(X,J)$ be a closed, reduced
almost complex $4$-orbifold. A $J$-holomorphic curve $C$ in $X$ is a
closed subset which is the image of a $J$-holomorphic map
$f:\Sigma\rightarrow X$. One can always
arrange that the following conditions hold.

\begin{itemize}
\item The map $f$ is not multiply covered, ie., the induced
map is not multiply covered.
\item For all but at most finitely many regular points $z\in\Sigma$, the
homomorphism $\rho_z$ in a local representative $(\hat{f}_z,\rho_z)$
of $f$ at $z$ is an
isomorphism. It is easily seen that with this assumption, the order
$m_\Sigma$ of the isotropy groups at regular points of $\Sigma$ depends
on $C$ only and equals the order of the isotropy group $G_p$ at a generic
point $p\in C$. We will denote this number by $m_C$ and call it the
`multiplicity' of $C$.
\end{itemize}
Such a $J$-holomorphic map $f$ is called a `parametrization' of $C$.
Note that $\Sigma$ may not be reduced, in which case
$m_\Sigma=m_C>1$. See \cite{C2} for more details.

Next we define an intersection product for any two $J$-holomorphic curves
and a pairing between a class in $H^2(X;\Q)$ and a $J$-holomorphic
curve, cf. \cite{C2}.

\begin{definition}
{\em (1)} For any $J$-holomorphic curve $C$ in $X$, the `Poincar\'{e}
dual' of $C$ is defined to be the class $PD(C)\in H^2(X;\Q)$ which is
uniquely determined by
$$
PD(C)\cup\alpha[X]=m_C^{-1}\alpha[C],\;\forall \alpha\in H^2(X;\Q),
$$
where $[C]$ is the class of $C$ in $H_2(X;\Z)$. {\em(}Note: $X$ and $C$ are
canonically oriented by $J$, and the cup product on $H^2(X;\Q)$ is
non-degenerate because $X$ is a $\Q$-homology manifold.{\em)}

{\em (2)} The intersection number of two $J$-holomorphic curves
$C,C^\prime$ {\em(}not necessarily distinct{\em)} is defined to be
$$
C\cdot C^\prime=PD(C)\cup PD(C^\prime)[X].
$$

{\em (3)} The pairing of a class $\alpha\in H^2(X;\Q)$ with a 
$J$-holomorphic curve $C$ is defined to be
$$
\alpha\cdot C=\alpha\cup PD(C)[X], \;\forall \alpha\in H^2(X;\Q).
$$
\end{definition}

The point here is that the fundamental class of an orbifold should be
counted only as a fraction of the fundamental class of the underlying
space when the orbifold is not reduced (cf. \cite{CR, C2}). Note that
Definition 2.1 coincides with the usual definition when $m_C=m_{C^\prime}=1$,
which is always the case when $X$ has only isolated singularities.

Now we digress on the local intersection number and local
self-intersection number of pseudoholomorphic discs in $\C^2$. To this end,
we recall the following local analytic property of pseudoholomorphic curves.
Let $D=\{z||z|<1\}$ be the unit disc in $\C$, and let
$f:(D,0)\rightarrow (\C^2,0)$ be a $J$-holomorphic map which is embedded
on $D\setminus\{0\}$. Then for any sufficiently small $\epsilon>0$, there
is an almost complex structure $J_\epsilon$ and a $J_\epsilon$-holomorphic
immersion $f_\epsilon$ such that as $\epsilon\rightarrow 0$,
$J_\epsilon\rightarrow J$ in $C^1$ topology and $f_\epsilon\rightarrow f$
in $C^2$ topology. Moreover, given any annuli $\{\lambda\leq |z|<1\}$
and $\{\lambda^\prime\leq |z|\leq\lambda\}$ in $D$, one can arrange
to have $f=f_\epsilon$ in $\{\lambda\leq |z|<1\}$ and to
have $J_\epsilon=J$ except in a chosen neighborhood of the image
of $\{\lambda^\prime\leq |z|\leq\lambda\}$ under $f$ by letting
$\epsilon>0$ sufficiently small. See \cite{McD1, McD2} for more details.

\begin{definition}{\em (See \cite{McD1, McD2}, and also \cite{C2})}
\begin{itemize}
\item [{(1)}] Let $C,C^\prime$ be any $J$-holomorphic discs in $\C^2$
parametrized by $f:(D,0)\rightarrow (\C^2,0)$, $f^\prime:(D,0)\rightarrow
(\C^2,0)$, such that both $f,f^\prime$ are embedded on $D\setminus
\{0\}$, and $C,C^\prime$ intersect at $0\in\C^2$ only. Then the local
intersection number $C\cdot C^\prime$ is defined using the following recipe:
perturb $f,f^\prime$ into $J_\epsilon$-holomorphic immersions
$f_\epsilon,f_\epsilon^\prime$, then
$$
C\cdot C^\prime=\sum_{\{(z,z^\prime)|f_\epsilon(z)
=f_\epsilon^\prime(z^\prime)\}}t_{(z,z^\prime)}
$$
where $t_{(z,z^\prime)}=1$ when $f_\epsilon(z)=f_\epsilon^\prime(z^\prime)$
is a transverse intersection, and $t_{(z,z^\prime)}=n\geq 2$ when
$f_\epsilon(z)=f_\epsilon^\prime(z^\prime)$ has tangency of order $n$. 
\item [{(2)}] Let $C$ be any $J$-holomorphic disc in $\C^2$ parametrized by
$f:(D,0)\rightarrow (\C^2,0)$ where $f$ is embedded on $D\setminus \{0\}$. Then
the local self-intersection number $C\cdot C$ is defined using the following
recipe: perturb $f$ into a $J_\epsilon$-holomorphic immersion $f_\epsilon$,
then
$$
C\cdot C=\sum_{\{[z,z^\prime]|z\neq z^\prime,
f_\epsilon(z)=f_\epsilon(z^\prime)\}}t_{[z,z^\prime]},
$$
where $[z,z^\prime]$ denotes the unordered pair of $z,z^\prime$, and
where $t_{[z,z^\prime]}=1$ when $f_\epsilon(z)=f_\epsilon(z^\prime)$
is a transverse intersection, and $t_{[z,z^\prime]}=n\geq 2$ when
$f_\epsilon(z)=f_\epsilon(z^\prime)$ has tangency of order $n$. 
\end{itemize}
\end{definition}

We remark that: (1) the local intersection number $C\cdot C^\prime$ 
depends only on the germs of $C,C^\prime$ at $0\in\C^2$, it is always 
positive, and $C\cdot C^\prime=1$ iff $C,C^\prime$ are both embedded 
and intersect at $0\in\C^2$ transversely, and (2) the local 
self-intersection number $C\cdot C$ depends only on the germ of $C$ 
at $0\in\C^2$, and it is non-negative with $C\cdot C=0$ iff $C$ is embedded.

Before stating the orbifold versions of the adjunction and intersection
formulae, we first introduce the relevant notation involved in the
statements.

(1) Let $C$ be a $J$-holomorphic curve parametrized by
$f:\Sigma\rightarrow X$. For any $z\in\Sigma$ with $p=f(z)$, let
$(\hat{f}_z,\rho_z)$ be a local representative of $f$ at $z$ and
$(V_p,G_p)$ be the local uniformizing system at $p$. Recall that
pseudoholomorphic curves (in manifolds) have only isolated
singular points (cf. \cite{McD1, McD2}). Since $f$ is not
multiply covered, it is easily seen that we may assume without
loss of generality that $\hat{f}_z$ is an embedding in
$V_p\setminus\{0\}$. (See \cite{C2} for more details.) With this
understood, we introduce
$$
\Lambda(C)_z\equiv \{\mbox{Im }(g\circ \hat{f}_z)|g\in G_p\},
$$
which is a set of $J$-holomorphic discs in $V_p\subset\C^2$
whose elements are naturally parametrized by the coset
$G_p/\mbox{Im }\rho_z$. We call an element of $\Lambda(C)_z$ a
`local representative' of $C$ at $z$. (In connection with $\Lambda(C)_z$,
observe that when $p$ is an isolated singular point, $(\hat{f}_z,\rho_z)$ can
be interpreted geometrically as follows: set $U_p\equiv V_p/G_p$, then
$\mbox{Im }\hat{f}_z\setminus\{0\}$ is a connected component of the inverse
image of $C\cap U_p\setminus\{p\}$ under the finite covering $V_p\setminus\{0\}
\rightarrow U_p\setminus\{p\}$, and $\mbox{Im }\rho_z$ is the subgroup of
deck transformations which leaves $\mbox{Im }\hat{f}_z$ invariant.)

(2) For any $J$-holomorphic curve $C$ in $X$, the `virtual genus' of $C$ is
defined to be
$$
g(C)=\frac{1}{2}(C\cdot C-c_1(TX)\cdot C)+\frac{1}{m_C}.
$$
Here $m_C$ is the multiplicity of $C$, ie., the order of $G_p$ at a
generic point $p\in C$.

(3) Let $\Sigma$ be an orbifold Riemann surface whose orbifold points have
orders $m_i$, $i=1,\cdots,k$, and whose regular points have order $m_\Sigma$,
and let $g_{|\Sigma|}$ be the genus of the underlying Riemann surface. Then
the `orbifold genus' of $\Sigma$ is defined to be
$$
g_\Sigma=\frac{g_{|\Sigma|}}{m_\Sigma}+\sum_{i=1}^k
(\frac{1}{2m_\Sigma}-\frac{1}{2m_i}).
$$
Note that with this definition, $c_1(T\Sigma)\cdot [\Sigma]=2m_{\Sigma}^{-1}
-2g_\Sigma$. (Here $[\Sigma]$ is $m_\Sigma^{-1}$ times the fundamental
class of the underlying Riemann surface, see the remarks below Definition 2.1.)

Now we state the adjunction and intersection formulae in \cite{C2}.

\begin{theorem}[Adjunction Formula]

Let $C$ be a $J$-holomorphic curve which is parametrized by
$f:\Sigma\rightarrow X$. Then
$$
g(C)=g_\Sigma + \sum_{\{[z,z^\prime]|z\neq z^\prime,
f(z)=f(z^\prime)\}}k_{[z,z^\prime]} +
\sum_{z\in\Sigma}k_z,
$$
where $[z,z^\prime]$ denotes the unordered pair of $z,z^\prime$, and
$k_{[z,z^\prime]},k_z$ are defined as follows.
\begin{itemize}
\item Let $G_{[z,z^\prime]}$ be the isotropy group at $f(z)=f(z^\prime)$
and $\Lambda(C)_z=\{C_{z,\alpha}\}$, $\Lambda(C)_{z^\prime}=
\{C_{z^\prime,\alpha^\prime}\}$, then
$$
k_{[z,z^\prime]}=\frac{1}{|G_{[z,z^\prime]}|}
\sum_{\alpha,\alpha^\prime}C_{z,\alpha}\cdot C_{z^\prime,\alpha^\prime}.
$$
\item Let $G_z$ be the isotropy group at $f(z)$ and $\Lambda(C)_z
=\{C_{z,\alpha}\}$, then
$$
k_z=\frac{1}{2|G_z|}(\sum_\alpha
C_{z,\alpha}\cdot C_{z,\alpha}+\sum_{\alpha,\beta}C_{z,\alpha}\cdot
C_{z,\beta}).
$$
\end{itemize}
\end{theorem}

(Note: the second sum $\sum_{\alpha,\beta}C_{z,\alpha}\cdot
C_{z,\beta}$ is over all $\alpha,\beta$ which are not necessarily distinct,
so it is the same as the self-intersection $(\sum_\alpha C_{z,\alpha})\cdot
(\sum_\alpha C_{z,\alpha})$.)

\begin{theorem}[Intersection Formula]

Let $C,C^\prime$ be distinct $J$-holomorphic curves parametrized by
$f:\Sigma\rightarrow X$, $f^\prime:\Sigma^\prime\rightarrow X$ respectively.
Then the intersection number
$$
C\cdot C^\prime=\sum_{\{(z,z^\prime)|f(z)=f^\prime(z^\prime)\}}
k_{(z,z^\prime)}
$$
where $k_{(z,z^\prime)}$ is defined as follows. Let $G_{(z,z^\prime)}$ be
the isotropy group at $f(z)=f^\prime(z^\prime)$ and $\Lambda(C)_z=
\{C_{z,\alpha}\}$, $\Lambda(C^\prime)_{z^\prime}=
\{C^\prime_{z^\prime,\alpha^\prime}\}$, then
$$
k_{(z,z^\prime)}=\frac{1}{|G_{(z,z^\prime)}|}
\sum_{\alpha,\alpha^\prime}C_{z,\alpha}\cdot
C^\prime_{z^\prime,\alpha^\prime}.
$$
\end{theorem}

(Note: the intersection formula is topological in nature in the sense
that it can be established for more general oriented surfaces, singular or
not, as long as the local intersection numbers
$C_{z,\alpha}\cdot C^\prime_{z^\prime,\alpha^\prime}$ can be
properly defined.)

Before we illustrate these formulae with examples, let's first
make some observations.

\begin{remark}

(1) The adjunction formula implies $g(C)\geq g_\Sigma$,
with $g(C)=g_\Sigma$ iff $k_{[z,z^\prime]}=0$ and $k_z=0$. The geometric
meaning of these conditions is that $C$ is topologically embedded and
at each singular point $p=f(z)$, $C$ has only one local representative
which is embedded (particularly, $\mbox{Im }\rho_z=G_z\equiv G_p$).
This is equivalent to saying that $C$ is a suborbifold.

(2) In many cases $\mbox{Im }\rho_z\neq G_p$, for instance, when $G_p$
is non-abelian and $p$ is isolated (hence $\mbox{Im }\rho_z$ is abelian).
It is useful to observe that in this case,
$$
k_z\geq \frac{1}{2m_z}(\frac{|G_p|}{m_z}-1)\geq \frac{1}{2m_z},
\mbox{ where $m_z$ is the order of the isotropy group at $z$.}
$$
Moreover, $k_z=\frac{1}{2m_z}(\frac{|G_p|}{m_z}-1)$ iff each local
representative of $C$ at $z$ is embedded and any two
different ones intersect transversely.

(3) Assume $X$ has only isolated singularities. Then a combination of the
adjunction and intersection formulae and the index formula may be used to
extract information about the local structure of a $J$-holomorphic curve
near a singular point. To be more precise, let
$(\hat{f}_z,\rho_z):(D_z,\Z_{m_z})\rightarrow (V_p,G_p)$ be a local
representative of $f$ at $z$ where $p=f(z)\in X$ is a singular point,
and let $0<m_1,m_2<m_z$ such that the action of $\mbox{Im }\rho_z$ on
$V_p$ is given by
$$
\rho_z(\mu_{m_z})\cdot (z_1,z_2)=(\mu_{m_z}^{m_1}z_1,\mu_{m_z}^{m_2}z_2),
\mbox{ where } \mu_m\equiv \exp(\sqrt{-1}\frac{2\pi}{m}).
$$
Then on the one hand, $m_1,m_2$ are constrained by the appearance in the
index formula associated to the dimension of the moduli space of
$J$-holomorphic curves (cf. the paragraph below Lemma 1.10), and on
the other hand, $m_1,m_2$ are related to $k_{[z,z^\prime]}$, $k_z$, and
$k_{(z,z^\prime)}$ as follows: By assuming $J$ integrable near the
singular points (note: this does not effect Lemma 1.10), we may write
$\hat{f}_z(w)=(u(w),v(w))$ for some holomorphic functions
$u(w)=aw^{l_1}+\cdots$ and $v(w)=bw^{l_2}+\cdots$, where
$l_i\equiv m_i\pmod{m_z}$, $i=1,2$. As shown in Lemma 2.6 below,
$k_{[z,z^\prime]}$, $k_z$, $k_{(z,z^\prime)}$ can be estimated in terms of
$l_1,l_2$, hence in terms of $m_1,m_2$.
\hfill $\Box$
\end{remark}

\begin{lemma}
\begin{itemize}
\item [{(1)}] Let $C$ be a holomorphic disc in $\C^2$ parametrized by
$f(z)=(a(z^{l_1}+\cdots), z^{l_2})$, which is embedded in
$\C^2\setminus\{0\}$. Then the local self-intersection
number satisfies $C\cdot C\geq\frac{1}{2}(l_1-1)(l_2-1)$.
\item [{(2)}] Let $C, C^\prime$ be distinct holomorphic discs in $\C^2$
which are parametrized by $f(z)=(a(z^{l_1}+\cdots), z^{l_2})$ and
$f^\prime(z)=(a^\prime(z^{l_1^\prime}+\cdots), z^{l_2^\prime})$
respectively, such that $f$, $f^\prime$ are embedded in $\C^2\setminus\{0\}$.
Then the local intersection number $C\cdot C^\prime\geq\min (l_1l^\prime_2,
l_2l^\prime_1)$. Here $l_1=\infty$ {\em(}resp. $l_1^\prime=\infty${\em)}
if $a=0$ {\em(}resp. $a^\prime=0${\em)}.
\end{itemize}
\end{lemma}

Lemma 2.6 is straightforward from the definition of local
intersection and self-intersection numbers in Definition 2.2, cf.
\cite{McD1, McD2}.

We end this section with the following example.

\begin{example}

Consider the weighted projective space $\P(d_1,d_2,d_3)$, which is
the quotient of $\C^3\setminus \{0\}$ under the $\C^\ast$-action
$$
z\cdot (z_1,z_2,z_3)=(z^{d_1}z_1,z^{d_2}z_2,z^{d_2}z_3), \;
\forall z\in\C^\ast\equiv \C\setminus\{0\}.
$$
Here we assume that $d_1,d_2,d_3$ are pairwise relatively prime and
satisfy $1<d_1<d_2<d_3$. Then $\P(d_1,d_2,d_3)$ is canonically a complex
$4$-orbifold, with three isolated singular points $p_1=[1:0:0]$,
$p_2=[0:1:0]$ and $p_3=[0:0:1]$ of isotropy group $\Z_{d_1}$, $\Z_{d_2}$
and $\Z_{d_3}$ respectively. Moreover, the set of orbifold complex
line bundles over $\P(d_1,d_2,d_3)$ is generated by the holomorphic
orbifold line bundle $E_0\equiv (\C^3\setminus \{0\})\times_{\C^\ast}\C$.
We will say that an orbifold complex line bundle $E$ (or the Poincar\'{e}
dual of $c_1(E)$) has degree $d$ if $E=E_0^{\otimes d}$. It is easy to check
that the canonical bundle $K$ has degree $-(d_1+d_2+d_3)$. Finally,
we observe that $c_1(E_0)\cdot c_1(E_0)=(d_1d_2d_3)^{-1}$.

With the preceding understood, let's apply the adjunction formula to
the holomorphic curve $C\equiv\{z_1=0\}$, which has degree $d_1$ and passes
through $p_2,p_3$. The virtual genus
$$
g(C)=\frac{1}{2}(C\cdot C+c_1(K)\cdot C)+1=\frac{1}{2}(\frac{d_1^2}{d_1d_2d_3}
+\frac{-(d_1+d_2+d_3)d_1}{d_1d_2d_3})+1=1-\frac{1}{2d_2}-\frac{1}{2d_3}.
$$
On the other hand, each singular point $p_2$ or $p_3$ will contribute
at least $\frac{1}{2}-\frac{1}{2d_2}$ or $\frac{1}{2}-\frac{1}{2d_3}$
to the right hand side of the adjunction formula. It follows easily
that $C$ is a suborbifold.

Next we look at the family of holomorphic curves
$$
C_\lambda\equiv \{az_2^{d_3}+bz_3^{d_2}=0\}, \mbox{ where }
\lambda=[a:b]\in\C\P^1\mbox{ with } ab\neq 0.
$$
Each $C_\lambda$ has degree $d_2d_3$ and passes through the singular point
$p_1$. To apply the adjunction formula to $C_\lambda$, we first calculate
$$
g(C_\lambda)=\frac{1}{2}(\frac{(d_2d_3)^2}{d_1d_2d_3}
-\frac{(d_1+d_2+d_3)d_2d_3}{d_1d_2d_3})+1=\frac{1}{2}-\frac{1}{2d_1}+
\frac{(d_2-1)(d_3-1)}{2d_1}.
$$
On the other hand, by Lemma 2.6 (1), the contribution $k_z$ at $p_1$
is at least $\frac{(d_2-1)(d_3-1)}{2d_1}$. It follows easily that
$C_\lambda$ is a $2$-sphere with only one singularity $p_1$.

Finally, we consider the intersection of $C_\lambda, C_{\lambda^\prime}$
where $\lambda\neq \lambda^\prime$. First, observe that $C_\lambda,
C_{\lambda^\prime}$ intersect at $p_1$ with local contribution
$k_{(z,z^\prime)}$ of at least $\frac{d_2d_3}{d_1}$ by Lemma 2.6 (2).
On the other hand, $C_\lambda\cdot C_{\lambda^\prime}
=\frac{(d_2d_3)^2}{d_1d_2d_3}=\frac{d_2d_3}{d_1}$. The intersection formula
implies that $C_\lambda\cap C_{\lambda^\prime}=\{p_1\}$.

All of the claims above can be easily verified directly.
\hfill $\Box$
\end{example}

\section{Seiberg-Witten theory and Taubes' theorems}

In this section, we first give a brief summary of the Seiberg-Witten theory
for smooth $4$-orbifolds, which is completely parallel to that for smooth
$4$-manifolds (see eg. \cite{T3, HT} for a nice introduction to Seiberg-Witten
theory of smooth $4$-manifolds). We then discuss the symplectic $4$-orbifold
analog of the theorems of Taubes in \cite{T1, T3} concerning nonvanishing of
the Seiberg-Witten invariant and the existence of certain pseudoholomorphic
curves. At the end of this section, we present an example which combines
various ingredients of the pseudoholomorphic curve theory and show
the existence of a symplectic $2$-sphere of degree $d_1$ or $-d_1$ in 
the weighted projective space
$\P(d_1,d_2,d_3)$ in Example 2.7 for any given symplectic structure.

All $4$-orbifolds in this section are assumed to be reduced, closed
and oriented. Let $X$ be a $4$-orbifold. We denote by $b_i(X)$ the dimension
of $H^i(X;\Q)$, and by $b_2^{+}(X)$ the dimension of the maximal
positive subspace for the cup product on $H^2(X;\Q)$.

Now let $X$ be a smooth 4-orbifold. Given any Riemannian metric on
$X$, a $Spin^\C$ structure is an orbifold principal $Spin^\C(4)$ bundle over
$X$ which descends to the orbifold principal $SO(4)$ bundle of oriented
orthonormal frames under the canonical homomorphism $Spin^\C(4)\rightarrow
SO(4)$. There are two associated orbifold $U(2)$ vector bundles (of rank $2$)
$S_{+},S_{-}$ with $\det(S_{+})=\det(S_{-})$, and a Clifford multiplication
which maps $T^\ast X$ into the skew adjoint endomorphisms of
$S_{+}\oplus S_{-}$.

The Seiberg-Witten equations associated to the $Spin^\C$ structure (if
there is one) are equations for a pair $(A,\psi)$, where $A$ is a connection
on $\det(S_{+})$ and $\psi$ is a section of $S_{+}$. Recall that the
Levi-Civita connection together with $A$ defines a covariant derivative
$\nabla_A$ on $S_{+}$. On the other hand, there are two maps $\sigma:
S_{+}\otimes T^\ast X\rightarrow S_{-}$ and $\tau:\mbox{End}(S_{+})
\rightarrow\Lambda_{+}\otimes\C$ induced by the Clifford multiplication,
with the latter being the adjoint of $c_{+}:\Lambda_{+}\rightarrow
\mbox{End}(S_{+})$, where $\Lambda_{+}$ is the orbifold bundle of
self-dual $2$-forms. With these understood, the Seiberg-Witten equations
read
$$
D_A\psi=0 \mbox{  and  } P_{+}F_A=\frac{1}{4}\tau(\psi\otimes\psi^\ast)
+\mu,
$$
where $D_A\equiv\sigma\circ\nabla_A$ is the Dirac operator, $P_{+}:
\Lambda^2 T^\ast X\rightarrow \Lambda_{+}$ is the orthogonal projection,
and $\mu$ is a fixed, imaginary valued, self-dual $2$-form  which
is added in as a perturbation term.

The Seiberg-Witten equations are invariant under the gauge transformations
$(A,\psi)\mapsto (A-2\varphi^{-1}d\varphi,\varphi\psi)$, where
$\varphi\in C^\infty(X;S^1)$ are circle-valued smooth functions on $X$.
The space of solutions modulo gauge equivalence, denoted by $M$, is compact,
and when $b^{+}_2(X)\geq 1$ and when it is nonempty, $M$ is a smooth
orientable manifold for a generic choice of $(g,\mu)$, where $g$ is the
Riemannian metric and $\mu$ is the self-dual $2$-form of perturbations.
Furthermore, $M$ contains no classes of reducible solutions (ie., those
with $\psi\equiv 0$), and if we let $M^0$ be the space of solutions modulo
the based gauge group, ie., those $\varphi\in C^\infty(X;S^1)$ such that
$\varphi(p_0)=1$ for a fixed base point $p_0\in X$, then $M^0\rightarrow M$
defines a principal $S^1$-bundle. Let $c$ be the first Chern class of
$M^0\rightarrow M$, $d=\dim M$, and fix an orientation of $M$. Then
the Seiberg-Witten invariant associated to the $Spin^\C$ structure is
defined as follows.
\begin{itemize}
\item When $d<0$ or $d=2n+1$, the Seiberg-Witten invariant is zero.
\item When $d=0$, the Seiberg-Witten invariant is a signed count of points
in $M$.
\item When $d=2n>0$, the Seiberg-Witten invariant equals $c^n[M]$.
\end{itemize}
As in the case of smooth 4-manifolds, the Seiberg-Witten invariant of $X$
is well-defined when $b^{+}_2(X)\geq 2$, depending only on the diffeomorphism
class of $X$ (as orbifolds). Moreover, there is an involution on the set of
$Spin^\C$ structures which preserves the Seiberg-Witten invariant up to
a change of sign. When $b^{+}_2(X)=1$, there is a chamber structure and the
Seiberg-Witten invariant also depends on the chamber which the pair $(g,\mu)$
is in. Moreover, the change of the Seiberg-Witten invariant when crossing
a wall of the chambers can be similarly analyzed as in the smooth 4-manifold
case.

Above is a brief summary of the Seiberg-Witten theory for smooth
$4$-orbifolds. In the case of global quotient $X=Y/G$, it corresponds to
the $G$-equivariant theory on $Y$.

Now we focus on the case where $X$ is a symplectic 4-orbifold. Let $\omega$
be a symplectic form on $X$. We orient $X$ by $\omega\wedge\omega$, and fix
an $\omega$-compatible almost complex structure $J$. Then with respect to
the associated Riemannian metric $g=\omega(\cdot,J\cdot)$, $\omega$ is
self-dual with $|\omega|=\sqrt{2}$. The set of $Spin^\C$ structures
on $X$ is nonempty. In fact, the almost complex structure $J$ gives rise to
a canonical $Spin^\C$ structure where the associated orbifold $U(2)$ bundles
are $S_{+}^0=\I\oplus K_X^{-1}$, $S_{-}^0=T^{0,1}X$. Here $\I$ is the trivial
orbifold complex line bundle and $K_X$ is the canonical bundle
$\det(T^{1,0}X)$. Moreover, the set of $Spin^\C$ structures is canonically
identified with the set of orbifold complex line bundles where each
orbifold complex line bundle $E$ corresponds to a $Spin^\C$ structure whose
associated orbifold $U(2)$ bundles are $S_{+}^E=E\oplus (K_X^{-1}\otimes E)$
and $S_{-}^E=T^{0,1}X\otimes E$. The involution on the set of $Spin^\C$
structures which preserves the Seiberg-Witten invariant up to a change of
sign sends $E$ to $K_X\otimes E^{-1}$.

As in the manifold case, there is a canonical (up to gauge equivalence)
connection $A_0$ on $K_X^{-1}=\det(S_{+}^0)$ such that the fact $d\omega=0$
implies that $D_{A_0}u_0=0$ for the section $u_0\equiv 1$ of $\I$ which is
considered as the section $(u_0,0)$ in $S_{+}^0=\I\oplus K_X^{-1}$.
Furthermore, by fixing such an $A_0$, any connection $A$ on
$\det(S_{+}^E)=K_X^{-1}\otimes E^2$ is canonically determined by a
connection $a$ on $E$. With these
understood, there is a distinguished family of the Seiberg-Witten equations
on $X$, which is parametrized by a real number $r>0$ and is for a triple
$(a,\alpha,\beta)$, where in the equations, the section $\psi$ of $S_{+}^E$ is
written as $\psi=\sqrt{r}(\alpha,\beta)$ and the perturbation term $\mu$
is taken to be $-\sqrt{-1}(4^{-1}r\omega)+P_{+}F_{A_0}$. (Here $\alpha$
is a section of $E$ and $\beta$ a section of $K_X^{-1}\otimes E$.)
Note that when $b_2^{+}(X)=1$, this distinguished family of Seiberg-Witten
equations belongs to a specific chamber for the Seiberg-Witten invariant.
This particular chamber will be referred to as the Taubes' chamber.
Now consider

\begin{theorem}
Let $(X,\omega)$ be a symplectic $4$-orbifold. Then the following are true.
\begin{itemize}
\item [{(1)}] The Seiberg-Witten invariant {\em (}in the Taubes' chamber
when $b_2^{+}(X)=1${\em)} associated to the canonical $Spin^\C$ structure
equals $\pm 1$. In particular, the Seiberg-Witten invariant corresponding
to the canonical bundle $K_X$ equals $\pm 1$ when $b_2^{+}(X)\geq 2$.
\item [{(2)}] Let $E$ be an orbifold complex line bundle. Suppose there
is an unbounded sequence of values for the parameter $r$ such that the
corresponding Seiberg-Witten equations have a solution $(a,\alpha,\beta)$.
Then for any $\omega$-compatible almost complex structure $J$, there are
$J$-holomorphic curves $C_1,C_2,\cdots,C_k$ in $X$ and positive integers
$n_1,n_2,\cdots,n_k$ such that $c_1(E)=\sum_{i=1}^k n_i PD(C_i)$. Moreover,
if a closed subset $\Omega\subset X$ is contained in $\alpha^{-1}(0)$
throughout, then $\Omega\subset\cup_{i=1}^k C_i$ also. 
\end{itemize}
\end{theorem}

(Here a $J$-holomorphic curve in $X$ is as defined in the beginning of 
\S 2. For the definition of $PD(C)$, the Poincar\'{e} dual of $C$, see 
Definition 2.1.)

Theorem 3.1 is the orbifold analog of the relevant theorems of Taubes
in \cite{T1, T3}. A proof of Theorem 3.1 can be found in \cite{C3}.

\begin{remark}

(1) Here is a typical source for the closed subset $\Omega$ in Theorem 3.1.
Suppose $p\in X$ is a singular point such that $G_p$ acts nontrivially on
the fiber of $E$ at $p$. Then $p\in\alpha^{-1}(0)$ for any solution
$(a,\alpha,\beta)$, and consequently $p\in\cup_{i=1}^k C_i$. Now let's
try this for the case where $X=Y/\Z_n$ and $E=K_X$. It follows that if
the space of $\Z_n$-equivariant self-dual harmonic $2$-forms on $Y$ has
dimension $\geq 2$, then for any given $\Z_n$-equivariant $J$ on $Y$,
$c_1(K_Y)=\sum_s n_s PD(C_s)$ for a finite set of
$J$-holomorphic curves $\{C_s\}$ in $Y$, such that
$\cup_s C_s$ is $\Z_n$-invariant, and contains any fixed point of
$\Z_n$ at which the complex representation of $\Z_n$ on the tangent space
has weights other than $(1,n-1)$.

(2) Even though, as we remarked in (1) above, that the $J$-holomorphic curves
$\{C_i\}$ contain certain singular points of $X$, nothing can be said
directly from the theorem about the local structure of these
curves near the singular points. To obtain such information, one may
use a combination of the adjunction and intersection formulae plus the
index formula, cf. Remark 2.5 (3). Here we point out another tactic,
which was used in \cite{C3}. Suppose the Seiberg-Witten invariant
corresponding to $E$ is nonzero and the dimension of the corresponding
Seiberg-Witten moduli space $M$ is $d=2n>0$. Then the $J$-holomorphic
curves $\{C_i\}$ may be required to contain any given subset $\Omega$
of less than or equal to $n$ points in $X$. One may further combine this
with the Gromov compactness theorem (cf. Theorem 1.13) by varying $\Omega$
to produce more desirable $J$-holomorphic curves.
\hfill $\Box$
\end{remark}

In light of Remark 3.2 (2), we give a formula for the dimension of the
Seiberg-Witten moduli space. For simplicity, we assume $X$ has only
isolated singular points. For the general case, a formula can be
found in Appendix A of \cite{C3}.

\begin{lemma}
Suppose the symplectic $4$-orbifold $X$ has only isolated singular points
$p_1,\cdots,p_l$. Then for any orbifold complex line bundle $E$, the dimension
$d(E)$ of the Seiberg-Witten moduli space corresponding to $E$ is given by
$$
d(E)=c_1(E)^2-c_1(E)\cdot c_1(K_X)+\sum_{i=1}^l I_i,
$$
where $K_X$ is the canonical bundle, and $I_i$ can be determined from the
singular point $p_i$ as follows. Let $\rho_{p_i,E}:G_{p_i}\rightarrow\C^\ast$
be the complex representation of $G_{p_i}$ on the fiber of $E$ at $p_i$, and
$\rho_{p_i,j}(g)$, $j=1,2$, be the eigenvalues of $g\in G_{p_i}$ associated
to the complex representation of $G_{p_i}$ on the tangent space at $p_i$. Then
$$
I_i=\frac{1}{|G_{p_i}|}\sum_{g\in G_{p_i}\setminus\{e\}}
\frac{2(\rho_{p_i,E}(g)-1)}{\prod_{j=1}^2(1-\rho_{p_i,j}(g^{-1}))},
\;\forall i=1,\cdots,l.
$$
\end{lemma}

We illustrate the formula with the following example.

\begin{example}

Let $E$ be the orbifold complex line bundle of degree
$d_1$ over the weighted projective space $\P(d_1,d_2,d_3)$ in
Example 2.7. We will show that $d(E)=0$.

Let $K$ be the canonical bundle. By Lemma 3.3,
$$
d(E)=c_1(E)^2-c_1(E)\cdot c_1(K)+I_1+I_2+I_3,
$$
where $I_1=0$, and with $\mu_m\equiv\exp(\sqrt{-1}\frac{2\pi}{m})$,
$$
I_2=\frac{1}{d_2}\sum_{x=1}^{d_2-1}\frac{2(\mu_{d_2}^{d_1x}-1)}
{(1-\mu_{d_2}^{-d_1x})
(1-\mu_{d_2}^{-d_3x})},\;
I_3=\frac{1}{d_3}\sum_{x=1}^{d_3-1}\frac{2(\mu_{d_3}^{d_1x}-1)}
{(1-\mu_{d_3}^{-d_1x})
(1-\mu_{d_3}^{-d_2x})}.
$$
We claim that
$$
I_2=\frac{d_2-1-2\delta_2}{d_2},\; I_3=\frac{d_3-1-2\delta_3}{d_3}
$$
where $\delta_2,\delta_3$ are uniquely determined by the
conditions $0\leq \delta_2\leq d_2-1$, $0\leq \delta_3\leq d_3-1$,
and $d_1-d_3\delta_2=0 \pmod{d_2}$, $d_1-d_2\delta_3=0\pmod{d_3}$.
Assuming validity of the claim, we have
\begin{eqnarray*}
d(E) & = &
\frac{d_1^2}{d_1d_2d_3}-\frac{-d_1(d_1+d_2+d_3)}{d_1d_2d_3}+
\frac{d_2-1-2\delta_2}{d_2}+\frac{d_3-1-2\delta_3}{d_3}\\
     & = & \frac{2(d_1+d_2d_3-\delta_3 d_2-\delta_2 d_3)}{d_2d_3}.
\end{eqnarray*}
Observe that (1) $\delta_3 d_2+\delta_2 d_3=d_1\pmod{d_2d_3}$ and
(2) $d_1<\delta_3 d_2+\delta_2 d_3< 2d_2d_3$, which implies that
$\delta_3 d_2+\delta_2 d_3=d_1+d_2d_3$, and hence $d(E)=0$.

It remains to verify the claim. We shall only check for
$I_2$, the case for $I_3$ is completely parallel. Introduce
$$
\varphi(t)\equiv\sum_{x=1}^{d_2-1}
\frac{\mu_{d_2}^{d_1x}}{1-\mu_{d_2}^{-d_3x}t}.
$$
Then by putting $\varphi(t)$ into power series of $t$, we have
\begin{eqnarray*}
\varphi(t) & = & \sum_{x=1}^{d_2-1}\mu_{d_2}^{d_1x}\sum_{l=0}^\infty
\mu_{d_2}^{-d_3xl}t^l\\
     & = & \sum_{l=0}^\infty\sum_{x=1}^{d_2-1}\mu_{d_2}^{(d_1-d_3l)x}t^l\\
     & = & \sum_{l=0}^\infty (-1)t^l+d_2t^{\delta_2}
\sum_{l=0}^\infty t^{d_2l}\\
     & = & (t-1)^{-1}+d_2t^{\delta_2}(1-t^{d_2})^{-1},
\end{eqnarray*}
which gives $\varphi(1)=\frac{1}{2}(d_2-1)-\delta_2$. Now it follows easily
that $I_2=\frac{2}{d_2}\varphi(1)=\frac{d_2-1-2\delta_2}{d_2}$.

\hfill $\Box$
\end{example}

We end this section by illustrating how the various ingredients are put into
use with the following example.

\begin{example}
We consider the $4$-orbifold $X=\P(d_1,d_2,d_3)$, the
weighted projective space in Example 2.7. Let $\omega$ be any
symplectic structure which evaluates positively on a homology
class in $H_2(X;\Q)$ of positive degree.

\vspace{2mm}

{\bf Claim 1:}\hspace{2mm}
{\em
For any $\omega$-compatible $J$, there exists a unique $J$-holomorphic
$2$-sphere of degree $d_1$, which is an embedded suborbifold containing
the singular points $p_2,p_3$ but not $p_1$.
}

\vspace{2mm}

The basic idea of the proof goes as follows. Let $E_1$ be the
orbifold complex line bundle of degree $d_1$. As in the smooth case, one
can show that the Seiberg-Witten invariant corresponding to $E_1$ is
nonzero in the Taubes' chamber, hence by Theorem 3.1, there exists a finite
set of $J$-holomorphic curves $\{C_\alpha\}$ such that
$$
c_1(E_1)=\sum_\alpha n_\alpha PD(C_\alpha) \mbox{ for some integers }
n_\alpha>0.
$$
We need to show that $\{C_\alpha\}$ consists of only one element
of degree $d_1$, from which Claim 1 follows easily with an
application of the adjunction formula. It is fairly easy to show that
$\{C_\alpha\}$ contains no more than two elements, but it is harder
to show that the single element has degree $d_1$. The point is that
the smaller the degree is, the larger the virtual genus is, and the
worse the situation is. The key observation here is that, on the other hand,
the more singular points the $J$-holomorphic curve contains, the larger
the right hand side of the adjunction formula is, and the better the
situation is. For this reason, we need to consider simultaneously the
orbifold complex line bundle $E_2$ of degree $d_2$.

Now the details of the proof. Let $K_\omega$ be the canonical bundle of
$(X,\omega)$. We will first show that $K_\omega$ has degree $-(d_1+d_2+d_3)$.
To see this, note that if $e(X), p_1(X)$ denote the Euler class and the first
Pontryagin class of $X$, then
$$
c_1(K_\omega)^2=2e(X)+p_1(X)=c_1(K)^2,
$$
where $K$ is the canonical bundle for the standard complex orbifold
structure on $X$, which has degree $-(d_1+d_2+d_3)$. It is clear that
$\deg K_\omega=-(d_1+d_2+d_3)$ if we show that $\deg K_\omega<0$, which
follows from $c_1(K_\omega)\cdot [\omega]<0$ because $[\omega]\cdot
c_1(E_1)>0$ by the assumption on $\omega$.
To see $c_1(K_\omega)\cdot [\omega]<0$, note that
$b_1(X)=0$ and $b_2^{+}(X)=b_2(X)=1$, so that in this case the chamber
structure for the Seiberg-Witten invariants is very simple:
there are only two chambers, the $0$-chamber and the Taubes' chamber.
The Seiberg-Witten invariant of the canonical $Spin^\C$ structure
in the Taubes' chamber is nonzero by Theorem 3.1, while the invariant is
zero in the $0$-chamber because $X$ has a metric of positive scalar
curvature (cf. \cite{BGN}.) Thus the wall-crossing number is nonzero,
which implies that $c_1(K_\omega)\cdot [\omega]<0$.

Next we claim that the Seiberg-Witten invariants corresponding to
$E_1, E_2$ are nonzero in the Taubes' chamber. To see this, note
that by the calculation in Example 3.4, the dimension $d(E_1)$
of the Seiberg-Witten moduli space is zero. A similar calculation shows
that $d(E_2)\geq 0$. It is easily seen that the claim follows from
a wall-crossing argument by observing that $c_1(E_i)\cdot [\omega]>0$,
$i=1,2$, and that $X$ has a metric of positive scalar curvature.

Now by Theorem 3.1, there exist finite sets of $J$-holomorphic curves
$\{C_\alpha\}$, $\{C_\beta\}$ such that
$$
c_1(E_1)=\sum_\alpha n_\alpha PD(C_\alpha),\;
c_1(E_2)=\sum_\beta n_\beta PD(C_\beta) \mbox{ for some integers }
n_\alpha, n_\beta>0.
$$
Moreover, it is easily seen that at $p_2 ,p_3$, the isotropy groups
$G_{p_2}$, $G_{p_3}$ act nontrivially on the fiber
of $E_1$, hence $p_2,p_3\in\cup_\alpha C_\alpha$ (cf. Remark 3.2 (1)).
Similarly, $p_1,p_3\in\cup_\beta C_\beta$.

Before we proceed further, let's make the following observation.

\vspace{2mm}

{\bf Claim 2:}\hspace{2mm}
{\em
The degree of a $J$-holomorphic curve in $X=\P(d_1,d_2,d_3)$ is always
integral.
}

\vspace{2mm}

Assuming the validity of Claim 2, we first consider the case where there
are $C_1\in\{C_\alpha\}$, $C_2\in\{C_\beta\}$ such that $C_1\neq C_2$.
Observe that $C_1\cap C_2\neq\emptyset$, because otherwise,
$C_1\cdot C_2=0$ which contradicts $b_2(X)=1$. Pick any point
$p\in C_1\cap C_2$, and
set $r_i\equiv \deg C_i$, $i=1,2$. Then by the intersection formula,
$$
\frac{r_1r_2}{d_1d_2d_3}=C_1\cdot C_2\geq\frac{1}{|G_p|}\geq\frac{1}{d_3},
$$
which implies that $r_1r_2\geq d_1d_2$. But $r_i\leq d_i$, $i=1,2$, so
that $r_1=d_1$, $r_2=d_2$. Now it is clear that in this case,
$\{C_\alpha\}$ consists of only $C_1$ with $\deg C_1=d_1$, which contains
both $p_2, p_3$. By the adjunction formula, $C_1$ is a $2$-sphere embedded
as a suborbifold in $X$ (see Example 2.7). Moreover, $C_1$ does not
contain the singular point $p_1$. Hence Claim 1 is proved in
this case except for the uniqueness. To see the uniqueness, suppose
there is a $C_1^\prime$ such that $\deg C_1^\prime=d_1$ and $C_1^\prime
\neq C_1$. Then similarly we have $C_1\cap C_1^\prime\neq\emptyset$, and
$$
\frac{d_1^2}{d_1d_2d_3}=C_1\cdot C_1^\prime\geq\frac{1}{d_3},
$$
which contradicts $d_1<d_2$. Hence the uniqueness.

There is still another possibility that both $\{C_\alpha\}$, $\{C_\beta\}$
consist of the same, single $J$-holomorphic curve $C_0$. But this possibility
can be ruled out as follows.

Suppose both $\{C_\alpha\}$, $\{C_\beta\}$ consist of $C_0$. Then $C_0$
contains all three singular points $p_1,p_2$ and $p_3$. Moreover,
$\deg C_0\in\Z$
by Claim 2, which implies that $\deg C_0=1$ because $\deg C_0$ is a
common divisor of $d_1,d_2$. Let $C_0$ be parametrized by a $J$-holomorphic
map $f:\Sigma\rightarrow X$. Since $p_1,p_2,p_3\in C_0$, there are $z_1,z_2,z_3
\in\Sigma$ such that $f(z_1)=p_1,f(z_2)=p_2$ and $f(z_3)=p_3$. Let
$m_1,m_2,m_3$ be the order of $z_1,z_2,z_3$ respectively.

First, we observe that if $m_i<d_i$ where $i=1,2$ or $3$, then the following
inequality holds for the contribution $k_{z_i}$ on the right hand side of the
adjunction formula:
$$
k_{z_i}\geq\frac{1}{2d_i}(\frac{d_i}{m_i}-1)\frac{d_i}{m_i}
\geq\frac{1}{2m_i}
$$
Second, we introduce the expression
$$
I\equiv 2(g(C_0)-1)+\frac{1}{d_1}+\frac{1}{d_2}+\frac{1}{d_3}.
$$
Then $\deg C_0=1$ implies that
$$
I=\frac{1}{d_1}+\frac{(d_2+d_3-1)(d_1-1)}{d_1d_2d_3}.
$$
It follows easily from the adjunction
formula that $I\geq 1$, which implies that
$$
\frac{1}{d_1}+\frac{d_2+d_3-1}{d_2d_3}> I\geq 1.
$$
But this contradicts the assumption that $d_1>1$ and
$d_1,d_2,d_3$ are pairwise relatively prime. Hence the second possibility
is ruled out.

It remains to verify Claim 2, which is a consequence of the adjunction
formula. More precisely, let $C$ be a $J$-holomorphic curve and let
$r=\deg C\in\Q$. Then the virtual genus
$$
g(C)=\frac{r^2-(d_1+d_2+d_3)r}{2d_1d_2d_3}+1.
$$
On the other hand, suppose $C$ is parametrized by a $J$-holomorphic
map $f:\Sigma\rightarrow X$ with orbifold genus ${g}_\Sigma=
g_{|\Sigma|}+\sum_{i=1}^k (\frac{1}{2}-\frac{1}{2m_i})$, where $m_i$
are the orders of the orbifold points on $\Sigma$. Then each $m_i$ is
a divisor of either $d_1,d_2$ or $d_3$. It is easily seen from the
adjunction formula that $r^2-(d_1+d_2+d_3)r\in\Z$, which implies that
$r=\deg C\in\Z$.

This completes the proof of Claim 1.

Let $C$ be the $J$-holomorphic $2$-sphere in Claim 1. It is interesting
to compare $C$ with the standard $2$-sphere $C_0\equiv\{z_1=0\}$.
In this regard, we have

\vspace{2mm}

{\bf Claim 3:}\hspace{2mm}
{\em $(X,C)$ and $(X,C_0)$ are diffeomorphic as orbifold pairs.
}

\vspace{2mm}

To see this, we first note that $C$ and $C_0$ have isomorphic normal bundles
hence have diffeomorphic regular neighborhoods in $X$. This is because the
bundles have the same Euler number (which is $C\cdot C=C_0\cdot C_0
=\frac{d_1}{d_2d_3}$) hence the same Seifert invariant, due to
the fact that $d_2,d_3$ are relatively prime. Secondly, the complements of the
regular neighborhoods are also diffeomorphic, ie., $C$ is not knotted in $X$.
This follows by observing that the complement $W$ of a regular neighborhood of
$C\sqcup \{p_1\}$ is a symplectic $\Z$-homology cobordism between lens spaces
(with a canonical contact structure on the boundary, see \S 4.1 below for more
details). By the main result in \cite{C2}, $W$ is smoothly a product, which
shows that $C$ is unknotted.
\hfill $\Box$
\end{example}

\section{Some applications}

This section is devoted to applications (or potential ones) of the
pseudoholomorphic curve theory for symplectic $4$-orbifolds.
More concretely, we will discuss some questions concerning the
following subjects where symplectic $4$-orbifolds are naturally
involved: (1) smooth $s$-cobordisms of elliptic $3$-manifolds,
(2) symplectic finite group actions on $4$-manifolds, (3) symplectic
circle actions on $6$-manifolds\footnote{I am grateful to Dusa McDuff
for communicating this problem to me.}, (4) algebraic surfaces with
isolated quotient singularities.

\subsection{Smooth $s$-cobordisms of elliptic $3$-manifolds}

An elliptic $3$-manifold is the quotient of $\s^3\subset\R^4$ under 
a free action
of some finite subgroup of $SO(4)$. An $s$-cobordism of elliptic
$3$-manifolds is a $4$-manifold with boundary the disjoint union
of two elliptic $3$-manifolds such that the inclusion of each
boundary component is a homotopy equivalence. (Note that in the
present case the Whitehead torsion vanishes, and the two boundary
components are homeomorphic, cf. \cite{KS1}.) An $s$-cobordism is
called trivial if it is the product of an elliptic $3$-manifold with
the closed interval $[0,1]$.

It has been a longstanding open question as whether there exist exotic
smooth structures on a trivial $s$-cobordism of elliptic $3$-manifolds
(cf. \cite{Ak}), or whether any of the nontrivial topological $s$-cobordisms
in \cite{CS, KS} is smoothable. In \cite{C3}, a ``geometrization'' program
was proposed to tackle these questions. At the center of the program is
the following conjecture.

\begin{conjecture}
A smooth $s$-cobordism of elliptic $3$-manifolds is smoothly trivial iff
its universal cover is smoothly trivial.
\end{conjecture}

The first step in the geometrization program, which was undertaken
in \cite{C3}, is to prove the triviality of a smooth $s$-cobordism of
elliptic $3$-manifolds under the assumption of existence of
a certain symplectic structure on the $s$-cobordism. More precisely,
given any $s$-cobordism of elliptic $3$-manifolds, which is
clearly orientable, one can fix an orientation and identify each
of the two boundary components with $\s^3/G\subset\C^2/G$ for some
finite subgroup $G\subset U(2)$ by an orientation-preserving
diffeomorphism (cf. \cite{C3}). Let $\omega_0$ be the canonical
symplectic structure on $\C^2/G$, which is the descendant of
$\sqrt{-1}(dz_1\wedge d\bar{z}_1+dz_2\wedge d\bar{z}_2)$ on $\C^2$.

\begin{theorem} {\em(See \cite{C3})}\hspace{2mm}
An $s$-cobordism of elliptic $3$-manifolds is smoothly trivial if
it has a symplectic structure which equals a positive multiple of
$\omega_0$ in a neighborhood of each boundary component.
\end{theorem}

(The case of lens space was treated in \cite{C2}, where it was
actually shown that a $\Z$-homology cobordism having such a symplectic
structure is smoothly trivial.)

The second step in the geometrization program is to show that a smooth
$s$-cobordism of elliptic $3$-manifolds must be symplectic if its
universal cover is smoothly trivial. It was explained in \cite{C3} that
this can be reduced to a special case of the problem of constructing
symplectic forms from generic self-dual harmonic ones discussed in Taubes
\cite{T4}. To be more precise, it was shown in \cite{C3} that
given any $s$-cobordism with each boundary component identified
with $\s^3/G\subset\C^2/G$ for some $G\subset U(2)$, there exists
a self-dual harmonic $2$-form which has regular zeroes and equals
a positive multiple of $\omega_0$ in a neighborhood of each
boundary component. It is easily seen that when the universal
cover is smoothly trivial, one can compactify the universal cover into
$\C\P^2$, and reduce the existence of a symplectic structure on the
$s$-cobordism to the following problem, cf. \cite{C3}.

\begin{problem}
Let $G$ be a finite group acting smoothly on $\C\P^2$ which has an
isolated fixed point $p$ and an invariant embedded $2$-sphere $S$
disjoint from $p$, such that $S$ is symplectic with respect to the
K\"{a}hler structure $\omega_0$ and generates the second homology.
Suppose $\omega$ is a $G$-equivariant, self-dual harmonic form which
vanishes transversely in the complement of $S$ and $p$. Modify $\omega$
in the sense of Taubes \cite{T4}, away from $S$ and $p$, to construct a
$G$-equivariant symplectic form on $\C\P^2$.
\end{problem}

Recall that Taubes' program in \cite{T4} is based on the following
two facts: (1) on an oriented $4$-manifold with $b_2^{+}\geq 1$, a
self-dual harmonic $2$-form has only regular zeroes for a generic
metric, which is a union of embedded circles, (2) if the
Seiberg-Witten invariant is nonzero, there are pseudoholomorphic
subvarieties in the complement of the zero set, which
homologically bound the zero set. Taubes in \cite{T4} asked under
what conditions one can modify the self-dual $2$-form in a neighborhood
of the pseudoholomorphic subvarieties to construct a symplectic
form. As explained in \cite{T4}, when the $4$-manifold is $S^1\times M^3$
and the self-dual $2$-form is $S^1$-equivariant, the problem is equivalent
to cancellation of all critical points of a circle-valued Morse function on
$M^3$ to make a fibration over $S^1$.

Although Taubes has since written a number of foundational papers
on this subject (cf. \cite{T5, T6, T7, T8,T9, T10}, see also
Gay-Kirby \cite{GK} and Auroux-Donaldson-Katzarkov \cite{ADK}),
the problem still remains largely open. However, if successful, it would
provide a very interesting alternative way of constructing symplectic
structures on $4$-manifolds, which is somewhat more abstract compared with
the existing methods of symplectic sum and Stein surface construction
(cf. \cite{Gom1, Gom2, McW0}). We hope that Conjecture 4.1 and
Problem 4.3 may serve as a testing ground for Taubes' problem as
well.

In connection with the last remark above, we mention the Cappell-Shaneson's
$s$-cobordism $W_{CS}$. Akbulut in \cite{Ak} showed that the universal cover
of $W_{CS}$ is smoothly trivial. Thus it would be very interesting to
determine by direct means (such as in \cite{Ak}) as whether $W_{CS}$ or any
of its two finite covers (other than the universal cover) is smoothly exotic.

A few words about the proof of Theorem 4.2. The case where the
elliptic $3$-manifold is $\s^3/G$ with $G$ nonabelian is more involved,
and its proof requires the full power of the pseudoholomorphic curve
theory. However, the basic idea of the proof can be simply described
as follows. Let $W$ be a symplectic $s$-cobordism of $\s^3/G$
where $G$ is nonabelian. We compactify $W$ into a symplectic
$4$-orbifold $X$ such that the concave end of $W$ is coned off and
the convex end is collapsed into a $2$-orbifold along the fibers
of the Seifert fibration. Then $X$ is homotopy equivalent to $X_0$,
where $X_0$ stands for $\B^4/G$ with boundary collapsed along the
Seifert fibration. It is shown that $X$ and $X_0$ are actually
diffeomorphic by recovering $X$ from the moduli space of pseudoholomorphic
curves in $X$ which corresponds to the family of complex lines in $\B^4/G$.
It follows easily from $X\cong X_0$ that $W$ is smoothly trivial.

\subsection{Symplectic finite group actions on $4$-manifolds}

Let $Y$ be a symplectic $4$-manifold and $G$ be a finite group
acting on $Y$ symplectically. Then the quotient space $Y/G$ is a
symplectic $4$-orbifold, and the pseudoholomorphic curve theory,
when applied to $Y/G$, may help extract information about the
action of $G$ on $Y$. The basic idea is that Taubes' existence theorem
will give a set of $J$-holomorphic curves in $Y$ which is invariant under
the group action, so that one may read off information about the group
action on the $4$-manifold $Y$ (particularly the fixed-point data)
from the induced action on the $2$-dimensional subset. 

One of the themes of symplectic topology is the comparison between
symplectic manifolds and K\"{a}hler manifolds. Extending to group
actions, it is natural to compare symplectic group actions on symplectic
manifolds with holomorphic actions on K\"{a}hler manifolds.

It is known that holomorphic actions on compact K\"{a}hler surfaces are
rigid in some aspects. For example, it was shown in \cite{Pe}
(compare also \cite{MN}) that for a large class of compact K\"{a}hler 
surfaces, a holomorphic action must
be trivial if it induces identity on the cohomology ring. (To the contrary,
it was demonstrated in \cite{E} that topological actions are extremely
flexible in this regard.) A study of homological rigidity of symplectic
finite group actions on $4$-manifolds was initiated in \cite{CK}.

For another rigidity phenomenon of holomorphic actions, we mention
a theorem of Gang Xiao \cite{Xiao} which generalizes a classical
result of Hurwitz on Riemann surfaces. Xiao's theorem says that
the order of automorphism group of a minimal algebraic surface of
general type is bounded by a multiple of the Chern number $c_1^2$.
It would be an interesting question as whether Xiao's theorem has
any symplectic analog.

Finally, we remark that a particular advantage of pseudoholomorphic 
curves over the holomorphic ones is that the almost complex structure
$J$ may be chosen generic. For an application of this aspect to
automorphisms of K\"{a}hler surfaces (in connection with Lemma 1.10 and 
Example 1.11), see the recent preprint \cite{C4}.

\subsection{Symplectic circle actions on $6$-manifolds}

It was shown in \cite{McD0} that to each symplectic $\s^1$-action on
a closed symplectic manifold $(X,\omega)$, one can associate a
generalized moment map $\mu:X\rightarrow\s^1$, which is defined
(up to a constant) by the equation $\omega(\xi,\cdot)=d\mu$, where $\xi=$
the vector field generating the $\s^1$-action. Moreover, one can analyze
the $\s^1$-action by looking at the reduced spaces
$\mu^{-1}(\lambda)/\s^1$, $\lambda\in\s^1$, which are symplectic orbifolds
in general for regular values of $\mu$. Particularly, the
$\s^1$-action is Hamiltonian if $\mu$ is not onto. When $X$ is
$6$-dimensional, the reduced spaces $\mu^{-1}(\lambda)/\s^1$ are
symplectic $4$-orbifolds for regular values $\lambda$, so the
pseudoholomorphic curve theory may be applied in this situation.

An Hamiltonian $\s^1$-action must have fixed points, which correspond to
the critical points of the moment map. A natural question is whether
the converse is true. It is known to be true for K\"{a}hler manifolds
(cf. \cite{Fr}).

The first counterexample was given by McDuff in \cite{McD0}, where the
author constructed a $6$-dimensional closed symplectic manifold
with a non-Hamiltonian $\s^1$-action which is not fixed-point free.
In the same paper, it was also proved that an $\s^1$-action on a
symplectic 4-manifold is Hamiltonian iff it has fixed points. The method
of the proof is to show that the generalized moment map is not onto,
which can be
verified by analyzing the change of the reduced spaces when crossing
a critical value of the generalized moment map. In general, the problem boils
down to finding an effective invariant which can tell the change of the
reduced spaces. In dimension $6$, it would be interesting to see if there is
any version of Gromov invariants of symplectic $4$-orbifolds that can be used
to detect the change.

In the aforementioned counterexample of McDuff, the fixed-point set is a
union of tori. It is a conceivable conjecture that an $\s^1$-action is
Hamiltonian if its fixed-point set is nonempty and isolated. So
far the conjecture has only been verified for semi-free actions
(cf. \cite{TW}).

\subsection{Algebraic surfaces with isolated quotient singularities}

Let $(X,\omega)$ be a symplectic $4$-orbifold with $b_2^{+}(X)\geq 2$.
By the orbifold version of Taubes' theorem (cf. Theorem 3.1 in \S 3),
the canonical class of $X$ is represented by a set of $J$-holomorphic
curves for any given $\omega$-compatible $J$. Moreover, a singular point
$p\in X$ is contained in these $J$-holomorphic curves as long as the
complex representation of the isotropy group $G_p$ on the tangent space is
not contained in $SL_2(\C)$ (cf. Remark 3.2 (1)). Thus the pseudoholomorphic
curve theory may be used to extract global restrictions on these singular
points of $X$.

An interesting class of symplectic $4$-orbifolds is provided by
algebraic surfaces with isolated quotient singularities. Indeed,
one can adopt the method in \cite{McW} of symplectic gluing along
contact hypersurfaces to show that an algebraic surface $X$ with isolated
quotient singularities has a natural orbifold symplectic structure such that
the corresponding canonical class is the same as that of $X$ as a complex
$4$-orbifold. Note that in this case, singularities $p\in X$ with $G_p$
not contained in $SL_2(\C)$ are precisely the non-Du Val singularities.

Algebraic surfaces with isolated quotient singularities form a very
natural class of singular spaces. (For a recent survey on this subject,
see \cite{KZ}.) In fact, quotient singularities in
dimension $2$ are nothing but the log terminal singularities, and the set of
algebraic surfaces with isolated quotient singularities form the
smallest category containing nonsingular surfaces, which is closed under
the log minimal model program (cf. \cite{Mat}).

\vspace{3mm}

\centerline{\bf Acknowledgements}

\vspace{3mm}

The author is grateful to S\l awomir Kwasik for many valuable
conversations and timely advice. He also wishes to extend his thanks 
to Dusa McDuff for many good suggestions after reading an earlier version 
of this article, and to an anonymous referee for his/her extensive comments 
which have led to much improved exposition.

\end{document}